\documentclass[12pt]{article}
\usepackage{amssymb,amsmath,amsthm,epsfig}
\usepackage{latexsym, enumerate}
\usepackage{eepic}
\usepackage{epic}
\usepackage{graphicx}
\usepackage{subfigure}
\usepackage{color}
\usepackage{ifpdf}
\usepackage{dsfont}
\usepackage{multirow}
\usepackage{array}
\usepackage{indentfirst}
\usepackage[ruled,vlined,linesnumbered]{algorithm2e}
\theoremstyle{plain}
\newtheorem{lem}{Lemma}[section]
\newtheorem{thm}[lem]{Theorem}

\theoremstyle{definition}

\theoremstyle{remark}
\newtheorem{rem}{Remark}[section]

\newcommand{\wt}{\widetilde}

\newcommand{\bb}{\mathbb}
\newcommand{\bm}{\boldsymbol}

\begin{document}
\title{ \large\bf A replica exchange preconditioned Crank-Nicolson Langevin dynamic MCMC method for Bayesian inverse problems}

\author{
Na Ou\thanks{School of Mathematics and Statistics, Changsha University of Science and Technology, Changsha, 410114, China (oyoungla@csust.edu.cn)}
\and
Zecheng Zhang\thanks{Department of Mathematics, Carnegie Mellon University, Pittsburgh, PA 15217, USA.}
\and
Guang Lin\thanks{Department of Mathematics $\&$ School of Mechanical Engineering, Purdue University, West Lafayete, IN 47907-2067, USA.({\tt guanglin@purdue.edu})}
}

\date{}
\maketitle
\begin{center}{\bf ABSTRACT}
\end{center}
The preconditioned Crank-Nicolson Markov Chain Monte Carlo (pCN-MCMC) method is one of the efficient approaches for large-scale Bayesian inference problems. It is constructed by discretizing the stochastic partial differential equation (SPDE); when the gradient information is coupled in the SPDE, the corresponding SPDE is also referred to as the Langevin diffusion. In this work, we present the replica exchange preconditioned Crank-Nicolson Langevin dynamic  (repCNLD) Monte Carlo method to speed up the convergence of the pCN-MCMC and tackle the challenge of multimodal distribution simulation. The replica exchange method comprises two preconditioned Langevin diffusion chains with low and high temperatures. \textcolor{black}{More specifically, the diffusion chain with high temperature aims at global exploration and the particles from the low-temperature diffusion are for the sake of local exploitation.}
\textcolor{black}{The particles from the two chains are allowed to exchange with a specific swapping rate}. We provide brief proof that the replica exchange preconditioned Langevin diffusion (repLD) has the same invariant distribution as the ones without swapping; the acceleration effect is also confirmed. We use the Crank-Nicolson scheme to discretize the repLD, and set the preconditioned matrix as the covariance matrix of the Gaussian prior. The discretization error is analyzed. Moreover, we use the multi-variance strategy to save the computational cost for the proposed method. \textcolor{black}{That is, we adopt a high-fidelity method to evaluate the forward model in the low-temperature chain while employing  a low-fidelity method to compute the forward model in the high-temperature chain.}
Therefore, perturbations are induced into the energy function, leading to the biased swapping rate.
\textcolor{black}{To correct the bias caused by the errors in the energy function, we provide an unbiased estimator derived from some assumptions on the forward solvers.} Furthermore, we provide the adjoint-based repCNLD algorithm under the discrete frame to speed up the simulations, where the gradient of the log-likelihood is calculated efficiently by the adjoint method. The proposed method is tested by sampling Gaussian mixture distribution and solving PDE-based nonlinear inverse problems numerically.


\section{Introduction}

Markov chain Monte Carlo (MCMC) approaches are widely used in many fields to study the statistical properties of complex systems and sample the posterior distribution in the Bayesian framework \cite{robert1999monte, liu2001monte}. The posterior usually involves high dimensionality, nonlinearity, and multimodality. The standard MCMC approaches, such as the random walk Metropolis-Hastings (RWMH), can mix arbitrarily slowly as the dimension of the unknowns increase \cite{beskos2009optimal, mattingly2012diffusion, cotter2013mcmc}, which means the number of forwarding model simulations may tend to infinity. Moreover, the nonlinearity of the map between the unknowns and output makes it computationally prohibitive for the Markov chain construction.

A family of dimension-independent MCMC algorithms was presented in \cite{beskos2008mcmc, cotter2013mcmc}, of which the proposal distribution was constructed by the preconditioned Crank-Nicolson (pCN) discretization of a stochastic partial differential equation (SPDE) that preserves the reference or target measure. When the proposal is generalized to the Langevin setting, the Metropolis-adapted Langevin algorithm (MALA) \cite{roberts1996exponential, robert2013monte} is resulted, which belongs to Langevin Monte Carlo (LMC) algorithms. Another branch of the LMC algorithms is the unadjusted Langevin algorithm (ULA)\cite{roberts1996exponential, durmus2019high, dalalyan2020sampling}, where the Markov chain is constructed by naturally discretizing the Langevin diffusion \cite{lemons1997paul, brooks2011handbook}. The ULA is computationally efficient for sampling high-dimensional probabilities due to the absence of the Metropolis-Hastings acceptance step \cite{2021Bayesian}. In this paper, we will use the ULA to build the Markov chains efficiently.

Note that the stochastic noise in the Langevin diffusion involved with the pCN-MCMC approach is scaled by a temperature parameter and a positive-definite preconditioned matrix; they greatly affect the efficiency of the ULA method; we focus on the effect of the temperature parameter in this paper. It plays a crucial role in global exploration and local exploitation trade. \textcolor{black}{When the temperature is low,
the particle can exploit the local structure of the likelihood landscape; however, it may trap in the local mode and hence fail to explore the global domain.
On the other hand, when the temperature is high, the particle can explore the global domain. However, it may not converge to the local optima.}

\textcolor{black}{To exploit the temperature effects and better use the advantages of the high and low temperatures, the authors in \cite{chen2020accelerating} propose the replica exchange Langevin diffusion (reLD) method.
It comprises two Langevin diffusion chains with two temperatures.
By allowing swapping the particle between two chains, the authors show that the low-temperature chain will converge to the same stationary distribution with fewer iteration steps.}

\textcolor{black}{Authors in \cite{deng2020non} later extend the reLD as an optimization algorithm to find the optima of the neural network.
In practice, similar to the gradient descent neural network optimization algorithms, it is impossible to send all samples to the GPU; the authors hence also assume a mini-batch setting.
Specifically, one can only obtain an estimation of the gradient with a stochastic batch of samples. As a result, they call their algorithms replica exchange stochastic gradient LD, or reSGLD.}

\textcolor{black}{It is known that the reLD can speed up the convergence, i.e., one may need fewer iterations to reach the stationary distribution.
However, since the algorithm employs two chains, the per-iteration cost doubles. Consequently, the total computation cost may be even higher than the single-chain algorithms.}
The multi-variance replica exchange method is developed in\cite{lin2021multi} to save computational cost. The method is then applied to handle noise data in the operator learning optimization \cite{lin2021accelerated}.
\textcolor{black}{The idea is to use two different estimations in two chains. Specifically, one can use a low-accuracy estimation in the high-temperature chain, which explores the region while adopting the normal estimator in the low-temperature chain. One hence can save costs in the high-temperature chain without sacrificing the convergence too much. }

We propose a replica exchange pCN Langevin dynamic (repCNLD) MCMC method. The replica exchange method comprises two preconditioned Langevin diffusion chains (repLD); these Langevin diffusion algorithms are assigned different temperature parameters. The particles from the two diffusion chains are allowed to swap with some specific probability. {\color{black}We prove that the repLD with swapping has the same invariant distribution with the ones without swapping, the acceleration effect of the repLD method is also confirmed in theory, the convergence of the diffusion chains is affected by the temperature parameter, preconditioned matrix, and the swapping intensity.}

To discretize the Langevin diffusion, we use the Crank-Nicolson algorithm \cite{cotter2013mcmc}. {\color{black}It is a linear-implicit scheme \cite{kloeden2001linear, beskos2008mcmc}, the implicitness parameter is incorporated only in the linear part of the drift to allow for analytical tractability. The Crank-Nicolson algorithm is the special case for the linear-implicit scheme introduced in \cite{beskos2008mcmc}. Since we set the preconditioned matrix to be the covariance matrix of the Gaussian prior, we refer to it as the preconditioned Crank-Nicolson scheme following \cite{cotter2013mcmc}. We analyze the discretization error; the mean squared error of discretization grows linearly against the step size.}

{\color{black}As the temperature parameter can be treated as an adjustment factor to the variance of the observed data given the predicted output, we adopt the multi-variance strategy \cite{lin2021multi} to save the doubled cost in two chains.} The forward model in the chain with low temperature is solved  with the high-fidelity method, while the chain with high temperature is assigned with a low-fidelity forward solver. This leads to the swapping rate bias, we correct the bias from the perspective of forwarding model approximation error. The modified factor in the swapping rate provides a guide for the choice of the low-fidelity model used in the second chain.

To further speed up the evaluation of the forward models, we provide the framework of the discrete adjoint method \cite{bardsley2018computational} to calculate the Jacobian of the output for nonlinear inverse problems. The derivative of the log-likelihood function is calculated through the forward model and the corresponding costate equation, which requires much less computational cost. The replica exchange methods are successfully used to escape from the local optima and explore distributions with multiple modes. The efficiency improvement and convergence acceleration effects are confirmed by some numerical examples. We summarize our contributions as follows.

\begin{enumerate}
    \item We propose a replica exchange preconditioned Langevin diffusion (repLD) method to handle the multimodal distribution simulation and speed up the sampling of the posterior density. The preconditioned matrix plays an important part in acceleration.
    \item We use the Crank-Nicolson algorithm to discrete the repLD, and set the preconditioned matrix as the covariance of the Gaussian prior. We also present the discretization error analysis.
    \item We propose an unbiased estimator of the swapping rate for the multi-variance repCNLD method, which provides a guide for the choice of the low-fidelity model used in the second chain.
    \item The framework of the discrete adjoint method is provided to efficiently calculate the gradient for PDE nonlinear inverse problems.
\end{enumerate}

The paper is organized as follows. In Section 2, we describe the setup of the Bayesian inference and the Crank-Nicolson scheme and present the repLDs and repCNLD in detail, some results are demonstrated. Section 3 considers the methods for simulation speeding up: the multi-variance strategy and the discrete adjoint method. Some numerical experiments are presented to demonstrate the efficiency of the proposed method in Section 4. Section 5 concludes with a summary discussion and tissues for future research.

\section{Replica exchange pCN Langevin dynamics MCMC}

\subsection{Bayesian inference}

Denote the unknown parameter as $\bm \xi \in {\bb R}^{n}$, let ${\bf G}: {\bm \xi}\rightarrow {\bb R}^{n_d}$ be the forward map, where $n_d$ is the dimension of the observed data. We assume that the unknown parameter and the measurement data ${\bf d}$ are related by
\begin{equation}\label{likelihood}
{\bf d}={\bf G} ({\bm \xi})+\varepsilon,
\end{equation}
where $\varepsilon \in {\bb R}^{n_d}$ is a vector of random noise and $\varepsilon \sim {\cal N}(0, \sigma_o^2I)$. Both ${\bf d}$ and $\bm \xi$ are random variables under the framework of Bayesian inference, {\color{black}as the Bayes' rule states,} the posterior distribution of $\bm \xi$ is
\begin{equation}\label{posterior}
\pi({\bm \xi}|{\bf d})\propto \pi({\bf d} | \bm \xi)\pi_{pri}(\bm \xi).
\end{equation}
The data enters the posterior through the likelihood function $\pi({\bf d} | \bm \xi)$, it is given by
\[
\pi({\bf d} | \bm \xi)= (2\pi\sigma_o^2)^{-\frac{n_d}{2}}\exp\bigg(-\frac{\|{\bf d}-{\bf G}(\bm \xi)\|_2^2}{2\sigma_o^2}\bigg),
\]
where $\| \cdot \|_2$ refers to the Euclidean norm. $\pi({\bf d} | \bm \xi)$ is the likelihood function and $\pi_{pri}(\bm \xi)$ is prior information of the unknown parameter $\bm \xi$. We focus on the Gaussian prior in this paper, which is widely used in statistical applications \cite{stein1999interpolation, hjort2010bayesian}. Without loss of generality, we assume the prior is ${\cal N}(m, B)$, i.e.,
\[
\pi_{pri}(\bm \xi)\propto \exp\bigg(-\frac{(\bm{\xi}-m)^TB^{-1}(\bm{\xi}-m)}{2}\bigg),
\]
where $m \in {\bb R}^{n}$ is the mean, and the positive and symmetric matrix $B\in \mathbb{R}^{n\times n}$ is the covariance matrix of the Gaussian prior.

Let
\[
\psi(\bm \xi)=-\log \pi({\bf d} | \bm \xi)
\]
be the potential and
\begin{equation}\label{energy}
U(\bm{\xi})=\frac{1}{2}(\bm{\xi}-m)^TB^{-1}(\bm{\xi}-m)+\psi(\bm{\xi})
\end{equation}
be the energy function, {\color{black}the posterior in the equation $(\ref{posterior})$ can then be expressed in the form}
\[
\pi({\bm \xi}| {\bf d})\propto \exp(-U(\bm{\xi})).
\]
{\color{black}From here on, we omit the conditioning on data ${\bf d}$ and write the posterior density as $\pi({\bm \xi})$.}

For some cases, when the forward map ${\bf G}$ is linear, the adoption of Gaussian priors lead to Gaussian posteriors \cite{franklin1970well}, the statistical estimators, e.g., the mean, the maximum a posteriori (MAP), covariance, credible intervals can be obtained easily. However, the forward map ${\bf G}$ is often nonlinear, resulting in the lack of an analytical expression of the posterior distribution, and poses a significant challenge to estimating the quantity of interests.

\subsection{The preconditioned Crank-Nicolson algorithms}
\label{pCN-A}

MCMC methods are often used to sample from the posterior distribution that arises from Bayesian inference. The preconditioned Crank-Nicolson MCMC (pCN-MCMC) method is one of the efficient approaches to dealing with large-scale inference problems. It is constructed by a discretization of stochastic partial differential equations (SPDEs)
\begin{equation}
\label{cn1}
{\rm d}\bm{\xi}_t=-L\left(B^{-1}(\bm{\xi}_t-m)+\alpha{\nabla \psi(\bm{\xi}_t)}\right){\rm d}t+\sqrt{2L\tau}{\rm d}\bm{W}_t,
\end{equation}
where $\nabla \psi$ is the negative derivative of the log-likelihood function, $L$ is an arbitrary symmetric, positive-definite matrix, ${\bm W}_t \in \mathbb{R}^n$ is a standard $n$-dimensional Brown motion and $\tau >0 $ is the temperature parameter. The SPDE  preserves the reference measure ${\cal N}(m, \tau B)$ for $\alpha=0$ or the target  measure
\[
\pi({\bm \xi})\propto \exp\left(-\frac{U(\bm{\xi})}{\tau}\right)
\]
for $\alpha=1$ \cite{stuart2004conditional, beskos2008mcmc, 2015Nonlinear}.

Our goal is to obtain the target measure through the natural discretization of the Langevin diffusion or SPDE. Hence we set $\alpha=1$ in the SPDE $(\ref{cn1})$. {\color{black} We use the liner-implicit scheme \cite{beskos2008mcmc} to obtain the discrete-time version of the continuous process, yielding,
\[
{\bm \xi}_{k+1}-{\bm \xi}_{k}=-\left((1-\theta)LB^{-1}{\bm \xi}_{k}+\theta LB^{-1}{\bm \xi}_{k+1}-LB^{-1}m+L\nabla \psi({\bm \xi}_{k})\right)\delta + \sqrt{2L\delta\tau}{\bm w}_k,
\]
for some $\theta \in (0, 1)$, ${\bm w}_k\sim {\cal N}(0,I)$, $I\in\mathbb{R}^{n\times n}$ is the identity matrix and $\delta$ is the time step. The formula can be rearranged as
\begin{equation*}
\left(I+\theta\delta LB^{-1}\right){\bm \xi}_{k+1}=\left(I-(1-\theta)\delta LB^{-1}\right){\bm \xi}_{k}+LB^{-1}\delta \left(m-B{\nabla \psi({\bm \xi}_{k})}\right)+\sqrt{2L\delta\tau}{\bm w}_k,
\end{equation*}
since the matrix $(I+\theta\delta LB^{-1})$ is symmetric positive-definite, the uniqueness of the location ${\bm \xi}_{k+1}$ is guaranteed given the realization of ${\bm \xi}_{k}$. Convergence of such implicit methods for SPDEs is investigated in \cite{1992Numerical}.

When the implicit parameter is set as $\theta=1/2$, we have the Crank-Nicolson scheme,
\begin{equation}\label{discrete}
\left(I+\frac{1}{2}\delta LB^{-1}\right){\bm \xi}_{k+1}=\left(I-\frac{1}{2}\delta LB^{-1}\right){\bm \xi}_{k}+LB^{-1}\delta\left(m-B{\nabla \psi({\bm \xi}_{k})}\right)+\sqrt{2L\delta\tau}{\bm w}_k.
\end{equation}
}
Moreover, we set the preconditioned matrix $L$ to be the covariance matrix $B$ of the Gaussian prior, the equation $(\ref{discrete})$ becomes
\begin{equation}\label{pCN-delta}
\left(1+\frac{1}{2} \delta\right){\bm \xi}_{k+1}=\left(1-\frac{1}{2}\delta\right){\bm \xi}_{k}+\delta\left(m-B\nabla \psi({\bm \xi}_{k})\right)+\sqrt{2B \delta\tau}{\bm w}_k.
\end{equation}
Let
\[
\beta=\frac{2\sqrt{2 \delta}}{2+ \delta},
\]
{\color{black}then for $\delta \in (0, 2)$, we have $\beta\in [0, 1]$.} The formula of the preconditioned Crank-Nicolson Langevin dynamics (pCNLD) is resulted,
\begin{equation}\label{pCN-langevin}
{\bm \xi}_{k+1}=\sqrt{1-\beta^2}{\bm \xi}_{k}+(1-\sqrt{1-\beta^2})(m-B\nabla \psi({\bm \xi}_{k}))+\beta \sqrt{B\tau}{\bm w}_k.
\end{equation}
Or equivalently, by replacing the term $B^{-1}({\bm \xi}_{k}-m)+\nabla \psi({\bm \xi}_{k})$ with $\nabla U({\bm \xi}_{k})$ in the equation $(\ref{pCN-delta})$, it can be expressed as
\begin{equation}\label{pCN-U}
{\bm \xi}_{k+1} = {\bm \xi}_{k}-(1-\sqrt{1-\beta^2}) B\nabla U\left({\bm \xi}_{k}\right)+\beta\sqrt{B\tau}\bm{w}_k.
\end{equation}
The pCNLD method belongs to the ULA, the MCMC samples generated by this method produce a biased approximation of the target distribution \cite{roberts1996exponential, dwivedi2018log}, where the bias is stemmed from the discretization error. The asymptotic bias can be reduced when the step size is sufficiently small.

\begin{rem}
For the Langevin dynamics discussed above, the Gaussian prior is incorporated into the posterior density; when the hybrid prior information is imposed, e.g., the Laplace and the Gaussian prior are jointly suggested, it is equivalent to set
\[
\psi(\bm{\xi})=-\log \pi(\bf d |\bm{\xi})+\gamma |\bm{\xi}|
\]
where $\gamma$ is the hyperparameter in the central Laplace distribution and $|\cdot|$ is the $l_1$ norm of the vector.
\end{rem}

\subsection{Replica exchange preconditioned Langevin diffusion}

We note that the temperature parameter $\tau$ and preconditioned matrix $L$ play crucial roles in the Langevin diffusion, we focus on the effect of the temperature parameter in this paper. \textcolor{black}{When the temperature is low, the particle $\bm{\xi}_t$ can exploit the local structure and converge locally fast, however, it may get trapped in the local optima and does not travel globally. As a result, when the solution landscape is multi-mode, the particle may not converge to the global optima in practice.} On the other hand,  when the temperature parameter $\tau$ is large, the convergence of the diffusion can be accelerated since the particles are able to explore the entire state space.


Inspired by the success of replica exchange Langevin diffusion presented in paper \cite{chen2020accelerating}, \textcolor{black}{we use two chains and allow swapping to tackle the multimodal problems and improve the pCNLD convergence efficiency.} Consider a pair of particles governed by two preconditioned Langevin diffusion chains as defined in $(\ref{cn1})$ with $\alpha=1$ and temperatures $\tau_1 < \tau_2$,
\begin{equation*}
\begin{split}
\mathrm{d}\bm{\xi}^{(1)}_t = -L\left(B^{-1}(\bm{\xi}^{(1)}_t-m)+\nabla \psi(\bm{\xi}^{(1)}_t)\right){\rm d}t+\sqrt{2L\tau_1}\rm{d}\bm{W}^{(1)}_t,\\
\mathrm{d}\bm{\xi}^{(2)}_t = -L\left(B^{-1}(\bm{\xi}^{(2)}_t-m)+\nabla \psi(\bm{\xi}^{(2)}_t)\right){\rm d}t+\sqrt{2L\tau_2}\rm{d}\bm{W}^{(2)}_t,
\end{split}
\end{equation*}
they can also be expressed equivalently as,
\begin{equation}
\label{cns-U}
\begin{split}
\mathrm{d}\bm{\xi}^{(1)}_t = -L\nabla U(\bm{\xi}^{(1)}_t){\rm d}t+\sqrt{2L\tau_1}{\rm d}\bm{W}^{(1)}_t,\\
\mathrm{d}\bm{\xi}^{(2)}_t = -L\nabla U(\bm{\xi}^{(2)}_t){\rm d}t+\sqrt{2L\tau_2}{\rm d}\bm{W}^{(2)}_t,
\end{split}
\end{equation}
where $\bm{W}^{(1)}$ is independent with $\bm{W}^{(2)}$. {\color{black}The equation $(\ref{cns-U})$ converges to the invariant distribution
\[
\mu(\bm{\xi}^{(1)}, \bm{\xi}^{(2)})= \pi(\bm{\xi}^{(1)}, \bm{\xi}^{(2)})\rm{d}\bm{\xi}^{(1)}\rm{d}\bm{\xi}^{(2)},
\]
with the density
\begin{equation}\label{invariant}
\pi(\bm{\xi}^{(1)}, \bm{\xi}^{(2)})\propto e^{-\frac{U(\bm{\xi}^{(1)})}{\tau_1}-\frac{U(\bm{\xi}^{(2)})}{\tau_2}}.
\end{equation}
When the two particles are allowed to swap from position $(\bm{\xi}_t^{(1)}, \bm{\xi}_t^{(2)})$ at time $t$ to position $(\bm{\xi}_{t+dt}^{(2)}, \bm{\xi}_{t+dt}^{(1)})$ with a swapping rate $cs(\bm{\xi}_t^{(1)}, \bm{\xi}_t^{(2)})dt$, where $c\geq 0$ is the swapping intensity, and
\begin{equation}\label{swap}
s(\bm{\xi}_t^{(1)}, \bm{\xi}_t^{(2)})=\min\bigg\{1, \exp\left(\tau_\delta\left(U(\bm{\xi}_t^{(1)})-U(\bm{\xi}_t^{(2)})\right)\right)\bigg\},
\end{equation}
where
$\tau_\delta=\frac{1}{\tau_1}-\frac{1}{\tau_2}$. We have the following result:}
\begin{thm}[Reversibility and invariance]\label{invariance}
For any $c$ with $cdt < 1$ satisfied, the infinitesimal generator of the Langevin jump process $\{\bm{\xi}_t\}_{t\geq 0}$ defined by equations $(\ref{cns-U})$ and $(\ref{swap})$ is
\[
\begin{array}{l}
\mathcal{L}^c\left(f(x, y)\right)=\underbrace{-[L\nabla_{x}U(x)] \cdot  \nabla_{x}f(x, y)+\tau_{1} L : \nabla_{x}\nabla_{x} f(x, y)}_{\mathcal{L}^c_1\left(f(x, y)\right)}\\
\underbrace{-[L\nabla_{y}U(y)] \cdot  \nabla_{y}f(x, y)+\tau_2 L : \nabla_{y}\nabla_{y} f(x, y)}_{\mathcal{L}^c_2\left(f(x, y)\right)}+\underbrace{ cs(x, y)\left(f(y, x)-f(x, y)  \right)}_{\mathcal{L}^c_s\left(f(x, y)\right)},
\end{array}
\]
the domain of the generator $\mathcal{D}(\mathcal{L}^c)$ is the space of all twice-differentiable functions with compact support. The symbol $\cdot$ and $:$ represent the Frobenius inner products, i.e., for any $x, y\in \mathbb{R}^n$,
\begin{eqnarray*}
  [L\nabla_{x}U(x)] \cdot  \nabla_{x}f(x, y) &=& \sum_i \left[L\nabla_{x}U(x)\right]_i \left(\nabla_{x}f(x, y)\right)_i \\
  L : \nabla_{x}\nabla_{x} f(x, y) &=& \sum_{i,j} \left[L\right]_{i,j}\frac{\partial^2 f(x, y)}{\partial x_i \partial x_j} .
\end{eqnarray*}
Furthermore, $\{\bm{\xi}_t\}_{t\geq 0}$ is reversible and its invariant distribution has density $(\ref{invariant})$.
\end{thm}
The proof of the Theorem \ref{invariance} is analog to the proof of the Lemma 3.2 in paper \cite{chen2020accelerating}, the term $\mathcal{L}_s^c$ arises from the swap has the same form as the one in \cite{chen2020accelerating}, the symmetry of the preconditioned matrix $L$ makes the equality formula involved with $\mathcal{L}_1^c$ and $\mathcal{L}_2^c$ established, a detailed proof is provided in the Appendix.

In the following, {\color{black} we will prove that in addition to the temperature parameters and the swapping, the preconditioned matrix also contributes to the acceleration of the replica exchange preconditioned Langevin diffusion chains (repLD).} Let $\nu_t$ be the distribution of the repLD $\{\bm{\xi}_t\}_{t\geq 0} $ at time $t$, $\mu$ be the invariant distribution. According to the analysis in \cite{chen2020accelerating}, the derivative of $\chi^2(\nu_t\| \mu)$ with time $t$ can be characterized by the Dirichlet form, i.e.,
\begin{equation}\label{divergence}
\frac{\mathrm{d}}{\mathrm{d} t} \chi^{2}\left(\nu_{t} \| \mu\right)=-2 \mathcal{E}^{c}\left(\frac{\mathrm{d} \nu_{t}}{\mathrm{d} \mu}\right),
\end{equation}
where the $\chi^2$-divergence is used to quantify the discrepancy between $\nu_t$ and $\mu$, it is defined as
\[
\chi^{2}\left(\nu_{t} \| \mu\right)=\int\left(\frac{\mathrm{d} \nu_{t}}{\mathrm{d} \mu}-1\right)^{2} \mathrm{d} \mu.
\]
The Dirichlet form is defined as
\[
\mathcal{E}^{c}(f)=\int \Gamma^{c}(f) \mathrm{d} \mu,
\]
where
\[
\Gamma^{c}(f)=\frac{1}{2}\left(\mathcal{L}^{c}\left(f^{2}\right)-2 f \mathcal{L}^{c}(f)\right)
\]
is the Carr\'{e} du Champ operator, $\mu$ and $\mathcal{L}^{c}$ is the invariant distribution and infinitesimal generator of $\{\bm{\xi}_t\}_{t\geq 0}$, respectively.
\begin{thm}[Acceleration]\label{acceleration}
For any fixed function $f$, $\mathcal{E}^{c}(f)$ has the form
\begin{equation}\label{Dirichlet}
\begin{aligned}
{\mathcal E}^{c}(f)&=\underbrace{ \int \left(\tau_{1}\left\|\sqrt{L}\nabla_{x} f(x, y)\right\|^{2}+\tau_{2}\left\|\sqrt{L}\nabla_{y} f(x, y)\right\|^{2}\right) \mathrm{d} \mu\left(x, y\right)}_{{\mathcal E}^0\geq 0} \\
&+\underbrace{\int \frac{c}{2} s\left(x, y\right)\left(f(y, x)-f(x, y)\right)^{2} \mathrm{d} \mu\left(x, y\right)}_{\geq 0},
\end{aligned}
\end{equation}
\end{thm}

The derivation of $\mathcal{E}^{c}(f)$ is provided in the Appendix. The first term $\mathcal{E}^0$ in $(\ref{Dirichlet})$ corresponds to the
repLD without swapping, larger temperatures $\tau_1$, $\tau_2$ and appropriate matrix $L$ may result in a larger $\mathcal{E}^0$, which means a faster rate of convergence. The nonnegativity of the second term in $(\ref{Dirichlet})$ can further promotes the Dirichlet form as analysed in \cite{chen2020accelerating} and therefore leads to faster evolution of $\chi^{2}\left(\nu_{t} \| \mu\right)$.

\subsection{Replica exchange pCN Langevin dynamics}

We discretize the repLD with the pCN scheme introduced in the equation $(\ref{pCN-langevin})$ or $(\ref{pCN-U})$, the pair of particles are then updated following
\begin{equation}\label{repCNLDs}
\begin{split}
{{\bm \xi}}^{(1)}_{k+1} = \sqrt{1-\beta^2}{{\bm \xi}}^{(1)}_{k}+\left(1-\sqrt{1-\beta^2}\right)\left(m-B\nabla \psi({{\bm \xi}}^{(1)}_k)\right)+\beta\sqrt{B\tau_1}\bm{w}^{(1)}_k,\\
{ {\bm \xi}}^{(2)}_{k+1} = \sqrt{1-\beta^2}{{\bm \xi}}^{(2)}_{k}+\left(1-\sqrt{1-\beta^2}\right)\left(m-B\nabla \psi({ {\bm \xi}}^{(2)}_k)\right)+\beta\sqrt{B\tau_2}\bm{w}^{(2)}_k,
\end{split}
\end{equation}
or
\begin{equation}\label{repCNLD-U}
\begin{split}
{\bm \xi}^{(1)}_{k+1} = {\bm \xi}^{(1)}_{k}-\left(1-\sqrt{1-\beta^2}\right) B\nabla U\left({\bm \xi}^{(1)}_{k}\right)+\beta\sqrt{B\tau_1}\bm{w}^{(1)}_k,\\
{ {\bm \xi}}^{(2)}_{k+1} = {{\bm \xi}}^{(2)}_{k}-\left(1-\sqrt{1-\beta^2}\right) B\nabla U\left({\bm \xi}^{(1)}_{k}\right)+\beta\sqrt{B\tau_2}\bm{w}^{(2)}_k,
\end{split}
\end{equation}
where the tunable parameter $\beta$ depends on the time step $\delta$, ${\bm w}^{(i)}_k\sim {\cal N}(0,I)$, for $i=1, 2$. We call the schemes above the replica exchange preconditioned Crank-Nicolson Langevin dynamics (repCNLD) method. In order to analyse the discretization error of the repCNLD method, {\color{black}we relax the terms $\delta$  and $\sqrt{\delta}$ from the $\beta$ involved parameter in the equation $(\ref{repCNLD-U})$, and obtain the equivalent form}
\begin{equation}\label{repCNLD-delta}
\begin{split}
{\bm \xi}^{(1)}_{k+1} = {\bm \xi}^{(1)}_{k}-\eta\delta  B\nabla U\left({\bm \xi}^{(1)}_{k}\right)+\eta\sqrt{2B\tau_1\delta}\bm{w}^{(1)}_k,\\
{ {\bm \xi}}^{(2)}_{k+1} = {{\bm \xi}}^{(2)}_{k}-\eta\delta B\nabla U\left({\bm \xi}^{(1)}_{k}\right)+\eta\sqrt{2B\tau_2\delta}\bm{w}^{(2)}_k,
\end{split}
\end{equation}
where $\eta=\frac{2}{2+\delta}$.

We have some assumptions for the energy function. Firstly, we assume that the energy function is $K$-smooth, that is, for all $x, y \in {\mathbb R}^d$, there exists a constant $K< +\infty$ such that
\begin{equation}\label{smooth}
\|\nabla U(x)-\nabla U(y)\| \leq K\|x-y\|,
\end{equation}
we also assume that there exist positive constants $\alpha_1$ and $\alpha_2$, such that for all $x\in {\mathbb R}^d$,
\begin{equation}\label{B-dissipative}
[B\nabla U(x)]\cdot x \geq \alpha_1\|x\|^{2}-\alpha_2.
\end{equation}

The repLD can be characterized as
\begin{equation}\label{continue}
\mathrm{d}\bm{\xi}_t = -H\nabla D(\bm{\xi}_t)\mathrm{d}t+\sqrt{H}\Sigma_t\mathrm{d}\bm{W}_t,
\end{equation}
where $\bm{\xi}_t=\left[\bm{\xi}_t^{(1)}; \bm{\xi}_t^{(2)}\right] \in {\mathbb R}^{2n}$, $D(\bm{\xi}_t)=\left[U(\bm{\xi}_t^{(1)}); U(\bm{\xi}_t^{(2)})\right]$, ${\bm W}_t \in {\mathbb R}^{2n}$ is the Brown motion, the matrix \( H=\left(\begin{array}{cc}
B & 0 \\
0 & B
\end{array}\right)\), \( \sqrt{H}=\left(\begin{array}{cc}
\sqrt{B} & 0 \\
0 & \sqrt{B}
\end{array}\right)\), here $B$ is the covariance matrix of the Gaussian prior, since the preconditioned matrix is set as $B$. $\Sigma_t$ is a random matrix which switches between \( \left(\begin{array}{cc}
\sqrt{2 \tau_{1} I} & 0 \\
0 & \sqrt{2 \tau_{2} I}
\end{array}\right)\)
and \(\left(\begin{array}{cc}
\sqrt{2 \tau_{2} I} & 0 \\
0 & \sqrt{2 \tau_{1} I}
\end{array}\right)\) with the probability $cs\left(\bm{\xi}_t^{(1)}, \bm{\xi}_t^{(2)}\right)\mathrm{d}t$.

Denote $\bm{\xi}^{\delta}(k)=\left[\bm{\xi}^{\delta, 1}(k); \bm{\xi}^{\delta, 2}(k)\right]$ as the solution of the discretization scheme $(\ref{repCNLDs})$ or $(\ref{repCNLD-U})$ at the $k$-th iteration, i.e., $\bm{\xi}^{\delta}(k)$ satisfies
\[
\bm{\xi}^{\delta}(k+1)={\bm \xi}^{\delta}(k)-\eta\delta  H\nabla D\left({\bm \xi}^{\delta}(k)\right)+\eta\sqrt{H\delta}\Sigma^\delta(k)\bm{w}(k),
\]
where $\bm{w}(k) \in {\mathbb R}^{2d}$ is a random vector generated from the standard Gaussian distribution. $\Sigma^\delta(k)$ is defined the same as $\Sigma_t$ but switches with the probability $cs\left(\bm{\xi}^{\delta, 1}(k), \bm{\xi}^{\delta, 2}(k)\right)\delta$. Let $\{\bm{\xi}^{\delta}_t\}_{t\geq 0}$ be the continuous-time interpolation of $\{\bm{\xi}^{\delta}(k)\}_{k\geq 0}$ which is defined as
\begin{equation}\label{continue-interpolation}
\bm{\xi}_t^{\delta}=\bm{\xi}_0-\eta\int_0^t H\nabla D\left({\bm \xi}^{\delta}_{\lfloor s/{\delta} \rfloor \delta}\right){\rm d}s+\eta \int_0^t  \sqrt{H}\Sigma^{\delta}_{\lfloor s/{\delta} \rfloor \delta} {\rm d}\bm{W}_s,
\end{equation}
here $\lfloor \cdot \rfloor$ denotes the floor function. It can be verified that given $\bm{\xi}_0={\bm \xi}^{\delta}(0)$, we then have
\[
\bm{\xi}_t^{\delta} = \bm{\xi}^{\delta}(k),
\]
for all integers $k>0$ and $t = k\delta$.
\begin{thm}[Discretization Error]\label{discrete-error}
When the energy function satisfies assumptions $(\ref{smooth})$ and $(\ref{B-dissipative})$, there exists a constant $\epsilon(n,\tau_1,\tau_2,c,K,\alpha_1,\alpha_2,T,B)$ that
depends on the dimension $n$, temperature parameters $\tau_1$ and $\tau_2$, swapping intensity $c$, smoothness
constant $K$ and constants $(\alpha_1, \alpha_2)$ in the inequality $(\ref{B-dissipative})$, length of the time interval $T$, and the covariance matrix of the Gaussian prior $B$, such that for all $t\in[0, T]$,
\[
{\bb E}\left[ \|\bm{\xi}_t-\bm{\xi}_t^\delta\|^2 \right]\leq \epsilon(n,\tau_1,\tau_2,c,K,\alpha_1,\alpha_2,T,B)\delta.
\]
provided the step size $\delta \in (0, 2)$ satisfies \( 0<\frac{2\delta}{2+\delta}<\frac{\alpha_1}{2\|B\|_F^2K^2}\).
\end{thm}
The proof of Theorem \ref{discrete-error} adapts from the proof framework of Theorem 3.6 in paper \cite{chen2020accelerating}, more details about the proof is demonstrated in the Appendix.

\section{Simulation speeding up}
Compared with the single-chain method, the simulation of the forwarding model for the replica exchange method doubles. Meanwhile, the gradient of the log-likelihood requires much computational resources. We provide two schemes to speed up the simulation of the proposed method: the multi-variance strategy and the discrete adjoint method.

\subsection{Multi-variance repCNLD}

\textcolor{black}{In practice, one only has an estimation of the energy function $U(\cdot)$ or its gradient.
In the optimization literature, this is due to the mini-batch setting. Specifically, the limited computing memory only allows one to compute the gradient by a random (stochastic) batch of samples. In the inverse PDE literature, errors caused by the forward solver and errors in the observation all contribute to the inaccurate estimation of the energy function \cite{2006precondition, chung2020multi}.
Consequently, investigators are interested in developing Langevin dynamic algorithms which assume a perturbation in the gradient estimation.
Following the convention in the optimization society, researchers usually call such kinds of algorithms stochastic gradient Langevin dynamic (SGLD).}
The authors of paper \cite{lin2021multi} proposed the multi-variance replica exchange stochastic gradient Langevin dynamics (m-reSGLD) method to accelerate the sampling process.

Similarly, for the proposed repCNLD method, we use the high-fidelity method to solve the forward model for the chain with low temperature, while in the high-temperature chain, the forward model is solved with the low-fidelity method. We refer to it as m-repCNLD method for the convenience of description. As discussed in paper \cite{deng2020non, lin2021multi}, the biased swapping rate could be modified when perturbations are induced into the energy function. Denote the energy function involved term in calculating the swapping rate as
\[
S(\bm{\xi}^{(1)}, \bm{\xi}^{(2)})=\exp\left(\tau_\delta\left(U(\bm{\xi}^{(1)})-U(\bm{\xi}^{(2)})\right)\right),
\]
it is then approximated by
\[
\begin{aligned}
\wt S({{\bm \xi}}^{(1)}, {{\bm \xi}}^{(2)})&=e^{\tau_\delta\left( U({ {\bm \xi}}^{(1)})-\wt U({ {\bm \xi}}^{(2)})\right)}\\
&=e^{\tau_\delta\left(U({ {\bm \xi}}^{(1)})-U({ {\bm \xi}}^{(2)})\right)}e^{\tau_\delta\left(U({ {\bm \xi}}^{(2)})-\wt U({ {\bm \xi}}^{(2)})\right)},\\
&=S({ {\bm \xi}}^{(1)}, { {\bm \xi}}^{(2)})e^{\tau_\delta\left(U({ {\bm \xi}}^{(2)})-\wt U({ {\bm \xi}}^{(2)})\right)},
\end{aligned}
\]
here $U(\cdot)$ and $\wt U(\cdot)$ are calculated by the accurate and approximated forward model, respectively. The expectation of $\wt S({{\bm \xi}}^{(1)}, { {\bm \xi}}^{(2)})$ is
\begin{equation}\label{bias}
\bb{E}\left[ \wt S({ {\bm \xi}}^{(1)}, { {\bm \xi}}^{(2)}) \right]=S({ {\bm \xi}}^{(1)}, { {\bm \xi}}^{(2)})\bb{E}\left[ e^{\tau_\delta\left(U({ {\bm \xi}}^{(2)})-\wt U({ {\bm \xi}}^{(2)})\right)} \right].
\end{equation}
Let $ {\bf \wt G}: {\bm \xi}\rightarrow {\bb R}^{n_d}$ be the approximated forward output for the second chain, ${\bf G}$ be the accurate forward output for the first chain, we assume the normality holds for the accurate and approximated model, i.e.,
\[
{\bf \wt G}({\bm \xi})\sim {\cal N}({\bf G}({\bm \xi}), \wt \sigma^2 I), \quad \text{for} \ \forall {\bm \xi}.
\]
The relationship between the observed data and the forward model output can be induced by equation $(\ref{likelihood})$,
\[
{\bf d}\sim {\cal N}({\bf G}({\bm \xi}) , \sigma_o^2I), \quad \text{for} \ \forall {\bm \xi}.
\]
We then have the difference between the accurate and approximated energy.
\[
\begin{aligned}
U({{\bm \xi}})-{\wt U}({{\bm \xi}})&=\psi({{\bm \xi}})-{\wt \psi}({{\bm \xi}})\\
&=\frac{1}{2\sigma_o^2}\left(\|{\bf d}-{\bf G}({\bm \xi})\|^2-\|{\bf d}-{\bf\wt G}({\bm \xi})\|^2\right)\\
&=\frac{1}{2\sigma_o^2}\left(\|{\bf G}({\bm \xi})-{\bf\wt G}({\bm \xi})\|^2-2\left({\bf d}-{\bf G}({\bm \xi})\right)^T\left({\bf G}({\bm \xi})-\bf{\wt G}({\bm \xi})\right) \right).
\end{aligned}
\]
{\color{black}The expectation in the right hand side of the equation $(\ref{bias})$ can be calculated as
\[
\bb{E}\left[ e^{\tau_\delta\left(U({ {\bm \xi}}^{(2)})-\wt U({ {\bm \xi}}^{(2)})\right)} \right]=\bb{E}\left[ e^{\frac{\tau_\delta}{2\sigma_o^2}\left(Z^TZ-2X^TZ\right)}\right],
\]
where
\[
Z={\bf G}({ {\bm \xi}^{(2)}})-{\bf\wt G}({ {\bm \xi}^{(2)}}), \quad
X={\bf d}-{\bf G}({ {\bm \xi}^{(2)}}),
\]
they are random variables independent with each other for any fixed ${\bm \xi}^{(2)}$, and $Z\sim {\cal N}(0, \wt \sigma^2 I)$, $X\sim {\cal N}(0, \sigma_o^2 I)$.
\begin{thm}\label{unbias}
Let $r_{\sigma}=\frac{\wt \sigma^2}{\sigma_o^2}$, then when $r_\sigma<\frac{1}{\tau_\delta^2+\tau_\delta}$ satisfied,
\begin{equation}\label{modiS}
\wt S_m({ {\bm \xi}}^{(1)}, {{\bm \xi}}^{(2)})=\left[1-(\tau_\delta+\tau_\delta^2)r_\sigma \right]^{\frac{n_d}{2}}e^{\tau_\delta( U({ {\bm \xi}}^{(1)})-\wt U({ {\bm \xi}}^{(2)}))},
\end{equation}
is an unbiased estimator of $S(\bm{\xi}^{(1)}, \bm{\xi}^{(2)})$.
\end{thm}
The Theorem \ref{unbias} can be easily proved, according to the definition of the expectation, we have
\[
\begin{aligned}
\bb{E}\left[ e^{\frac{\tau_\delta}{2\sigma_o^2}\left(Z^TZ-2X^TZ\right)}\right]&=\int  \int (2\pi \wt \sigma^2)^{-\frac{n_d}{2}} e^{-\frac{\sigma_o^2-\wt \sigma^2(\tau_\delta^2+ \tau_\delta)}{2\sigma_o^2\wt \sigma^2}z^Tz}
(2\pi \sigma_o^2)^{-\frac{n_d}{2}} e^{-\frac{(x+\tau_\delta z)^T(x+\tau_\delta z)}{2\sigma_o^2}}{\rm d}x{\rm d}z\\
&=\int (2\pi \wt \sigma^2)^{-\frac{n_d}{2}} e^{-\frac{\sigma_o^2-\wt \sigma^2(\tau_\delta^2+ \tau_\delta)}{2\sigma_o^2\wt \sigma^2}z^Tz}{\rm d}z\\
&=\left[1-(\tau_\delta+\tau_\delta^2)r_\sigma \right]^{-\frac{n_d}{2}},
\end{aligned}
\]
which indicates that the term $\wt S({{\bm \xi}}^{(1)}, {{\bm \xi}}^{(2)})$ is a biased estimator of $S({ {\bm \xi}}^{(1)}, { {\bm \xi}}^{(2)})$. It can be verified that
\[
\begin{aligned}
\bb{E}\left[\wt S_m({ {\bm \xi}}^{(1)}, {{\bm \xi}}^{(2)})\right]&= \left[1-(\tau_\delta+\tau_\delta^2)r_\sigma \right]^{\frac{n_d}{2}}\bb{E}\left[ e^{\tau_\delta( U({ {\bm \xi}}^{(1)})-\wt U({ {\bm \xi}}^{(2)}))}\right]\\
&=\left[1-(\tau_\delta+\tau_\delta^2)r_\sigma \right]^{\frac{n_d}{2}}\bb{E}\left[ \wt S({ {\bm \xi}}^{(1)}, { {\bm \xi}}^{(2)}) \right]\\
&=S({ {\bm \xi}}^{(1)}, { {\bm \xi}}^{(2)}).\\
\end{aligned}
\]
 }
When applying the m-repCNLD method, the term $\wt S_m({ {\bm \xi}}^{(1)}, {{\bm \xi}}^{(2)})$ would be used to calculate the swapping rate.

\subsection{Discrete adjoint method for gradient calculation}

We note that the map between the unknowns ${\bm \xi}$ and output of the PDEs is nonlinear, the gradient of the log-likelihood function with respect to the unknown parameters is required, we use the adjoint method \cite{bardsley2018computational} to calculate the gradient efficiently. Let the output of the forward model be expressed in the following form.
\[
{\bf G}({\bm \xi})={\bf C A}({\bm \xi})^{-1}{\bf F}({\bm \xi}),
\]
where ${\bf C} \in {\bb R}^{n_d\times N}$ defines the linear state-to-observation map and ${\bf A}({\bm \xi})^{-1}{\bf F}({\bm \xi})$ is the solution of the discretized PDE, where ${\bf A}({\bm \xi}) \in {\bb R}^{N\times N}$ is assumed to be invertible and Fr\'{e}chet differentiable \cite{bardsley2018computational}, ${\bf F}({\bm \xi})$ is also assumed to be Fr\'{e}chet differentiable here.

Denote ${\bf r}(\bm \xi)={\bf G}({\bm \xi})-{\bf d}$ and ${\bf J}(\bm \xi)$ the Jacobian of ${\bf G}({\bm \xi})$, {\color{black}recall that
\[
\psi(\bm \xi)=-\log \pi({\bf d} | \bm \xi)=\frac{n_d}{2}\log 2\pi\sigma_o^2+\frac{\|{\bf d}-{\bf G}(\bm \xi)\|_2^2}{2\sigma_o^2},
\]
we then have }
\begin{equation}\label{gradient}
\nabla \psi(\bm \xi)=\frac{1}{\sigma_o^2}{\bf J}(\bm \xi)^T{\bf r}(\bm \xi),
\end{equation}
notice that the $i$th row of ${\bf J}(\bm \xi)^T{\bf r}(\bm \xi)$ can be expressed as
\begin{equation}\label{irow}
\begin{aligned}
{[{\bf J}(\bm \xi)^T{\bf r}(\bm \xi) ]_i}&=\frac{\partial {\bf G}({\bm \xi})}{\partial \xi_i} \cdot {\bf r}(\bm \xi) \\
&=-\left[{\bf C A}({\bm \xi})^{-1}\left({\bf F}_i^{\prime}({\bm \xi})- {\bf A}_i^{\prime}({\bm \xi}){\bf A}({\bm \xi})^{-1}{\bf F}({\bm \xi}) \right)\right] \cdot {\bf r}(\bm \xi),
\end{aligned}
\end{equation}
where the symbol $\cdot$ represents the Frobenius inner products, ${\bf A}_i^{\prime}({\bm \xi}):=\frac{d}{d\tau}{\bf A}({\bm \xi}+\tau {\bf e_i})|_{\tau=0}$, ${\bf F}_i^{\prime}({\bm \xi}):=\frac{d}{d\tau}{\bf F}({\bm \xi}+\tau {\bf e_i})|_{\tau=0}$. {\color{black}The second equality in the equation $(\ref{irow})$ holds because
\begin{equation}\label{GD}
\frac{\partial {\bf G}({\bm \xi})}{\partial \xi_i}={\bf C}\frac{\partial {\bf A}^{-1}({\bm \xi})}{\partial \xi_i}{\bf F}({\bm \xi})+{\bf C}{\bf A}^{-1}({\bm \xi}){\bf F}_i^{\prime}({\bm \xi}),
\end{equation}
since
\[
{\bf A}({\bm \xi}){\bf A}^{-1}({\bm \xi})=I,
\]
we have
\[
\frac{\partial {\bf A}^{-1}({\bm \xi})}{\partial \xi_i}=-{\bf A}^{-1}({\bm \xi}){\bf A}_i^{\prime}({\bm \xi}){\bf A}^{-1}({\bm \xi}),
\]
substitute it into the equation $(\ref{GD})$, the corresponding expression of  $\frac{\partial {\bf G}({\bm \xi})}{\partial \xi_i}$ is resulted.
}

Let ${\bf v}$ and ${\bf z}$ satisfy the forward model
\begin{equation}\label{Forward}
{\bf A}({\bm \xi}){\bf v}={\bf F}({\bm \xi}),
\end{equation}
and the costate equation
\begin{equation}\label{cos}
{\bf A}({\bm \xi})^T{\bf z}=-{\bf C}^T{\bf r}(\bm \xi),
\end{equation}
respectively. The equation $(\ref{irow})$ can then be arranged as
\[
{[{\bf J}(\bm \xi)^T{\bf r}(\bm \xi) ]_i}=\left[{\bf A}_i^{\prime}({\bm \xi}){\bf v}-{\bf F}_i^{\prime}({\bm \xi})\right] \cdot {\bf z},
\]
it means that for any ${\bm \xi}$, we need only solve the forward model and the costate equation to compute $\nabla \psi(\bm \xi)$. Meanwhile, the matrices ${\bf A}_i^{\prime}({\bm \xi})$ and ${\bf F}_i^{\prime}({\bm \xi})$ can be calculated in parallel. We list the discrete adjoint-based repCNLD method in Algorithm $\ref{repCNLD}$, the {\em for-loop} can be proceeded in parallel to save the computational cost. For the case of m-repCNLD, the forward model and costate equations are solved with the low-fidelity method for the second chain, the term in the swapping rate should be modified as equation $(\ref{modiS})$.

\begin{algorithm}
	\caption{Discrete adjoint-based repCNLD method.}
	\label{repCNLD}
	\KwIn{The mean $m$ and covariance matrix $B$ of the prior, parameter $\beta$, $\tau_1$, $\tau_2$, observed data ${\bf d}$, observation matrix $\bf C$, noise of the observed data $\sigma_o$. The number of iterations $N_{max}$. {\color{black}The intensity $c$ and time step $\delta$ are omitted in the swapping rate.} }
	\KwOut{Samples of Markov chains $\{\Xi_i^{(1)}\}_{i=1}^{N_{max}}$ and $\{\Xi_i^{(2)}\}_{i=1}^{N_{max}}$}
	\BlankLine
	Generate the initial values $\xi_0^{(1)}$, $\xi_0^{(2)}$ for the two chains randomly.

  	Evaluate the forward model at $\xi_0^{(j)}$ by equation $(\ref{Forward})$, obtain the solution ${\bf v}_0^{(j)}$ and errors $ {\bf r}_0^{(j)}={\bf C}{\bf v}_0^{(j)}-{\bf d} $, for $j=1, 2$.

	\While{$i<N_{\max}$}{	
		\For {$j=1:2$}{

         Solve the costate equation $(\ref{cos})$ and calculate the gradient term $\nabla \psi(\xi_0^{(j)})$ by equation $(\ref{gradient})$.

         Update the new particles by
           \[
           \xi^{(j)} = \sqrt{1-\beta^2}\xi_0^{(j)}+(1-\sqrt{1-\beta^2})(m-B\nabla \psi(\xi_0^{(j)}))+\beta\sqrt{\tau_j}w^{(j)}, \ w^{(j)}\sim \mathcal{N}(0, B).
           \]

       Simulate the forward model at $\xi^{(j)}$ by equation $(\ref{Forward})$, obtain the solution ${\bf v}^{(j)}$ and error $ {\bf r}^{(j)}={\bf C}{\bf v}^{(j)}-{\bf d} $.

       Calculate the energy function $E^{(j)}=U(\xi^{(j)})$ according to equation $(\ref{energy})$.

       }

       Draw $u$ from the uniform distribution $\mathcal{U}[0, 1]$ and calculate the exchange rate
        \[
        s=\min\bigg\{1, \exp\left(\tau_\delta(E^{(1)}-E^{(2)})\right)\bigg\},
        \]

        \If {$u< s$} {
              Update $\xi_0^{(1)}$ by $\xi^{(2)}$, $\xi_0^{(2)}$ by $\xi^{(1)}$.

              Update ${\bf v}_0^{(1)}$ by ${\bf v}^{(2)}$, ${\bf v}_0^{(2)}$ by ${\bf v}^{(1)}$.

              Update ${\bf r}_0^{(1)}$ by ${\bf r}^{(2)}$, ${\bf r}_0^{(2)}$ by ${\bf r}^{(1)}$.

              {\bf else}{

              Update $\xi_0^{(j)}$ by $\xi^{(j)}$, for $j=1,2$.

              Update ${\bf v}_0^{(j)}$ by ${\bf v}^{(j)}$, for $j=1,2$.

              Update ${\bf r}_0^{(j)}$ by ${\bf r}^{(j)}$, for $j=1,2$.
              }
              }
        Save the samples as \(
        \Xi_i^{(1)}=\xi_0^{(1)}, \ \ \Xi_i^{(2)}=\xi_0^{(2)}.
        \)
	}
\end{algorithm}

\section{Numerical experiments}
In this section, \textcolor{black}{we present some numerical experiments results. In particular, }we apply the proposed repCNLD method to sample Gaussian mixture distribution and solve Bayesian nonlinear inverse problems. \textcolor{black}{The Gaussian mixture distribution simulation shows that the method can capture multi-mode but not get stuck in local modes.} For the nonlinear PDE inverse problems, the adjoint method is used to compute the gradient of the logarithmic likelihood function effectively. \textcolor{black}{The proposed methods are able to explore the multi-modal and high-dimensional posterior density,} meanwhile improving the efficiency of the MCMC method.

%

\subsection{Gaussian mixture density simulation}
In this subsection, we will show the performance of the proposed method in sampling the Gaussian mixture density, i.e., for ${\bm \xi} \in {\bb R}^n$, the target distribution is
\[
\pi({\bm \xi}|\theta)=\sum_{k=1}^{M} \gamma_k\phi({\bm \xi}|\theta_k),
\]
where $\gamma_k \geq 0$ is the weight coefficients and $\sum_{k=1}^{M} \gamma_k = 1$, $\phi({\bm \xi}|\theta_k)$ is the Gaussian density and has the form
\[
\phi({\bm \xi}|\theta_k)=(2\pi)^{-\frac{n}{2}}|B_k|^{-\frac{n}{2}}\exp\bigg\{-\frac{1}{2}({\bm \xi}-m_k)^TB_k^{-1}({\bm \xi}-m_k)\bigg\},
\]
where $\theta_k = (m_k, B_k)$ is the given parameters of the density.
{\color{black}We use the scheme $(\ref{pCN-langevin})$ and $(\ref{repCNLDs})$ for the pCNLD and repCNLD method, respectively. The potential function is set as
\begin{equation}\label{m-potential}
\psi({\bm \xi})=-\log \pi({\bm \xi}|\theta)-\frac{1}{2}(\bm{\xi}-m)^TB^{-1}(\bm{\xi}-m),
\end{equation}
such that $\exp\{-U({\bm \xi})\}$ equals to the target distribution $\pi({\bm \xi}|\theta)$. Substitute the equation $(\ref{m-potential})$ into the equation $(\ref{pCN-langevin})$, the updating scheme for a single chain can be rearranged as
\[
{\bm \xi}_{k+1} ={\bm \xi}_{k}+(1-\sqrt{1-\beta^2})B\nabla \log \pi({\bm \xi}_k|\theta)+\beta\sqrt{B\tau}{\bm w}_k,
\]
the scheme for the repCNLD method can be correspondingly derived. It can be seen that the schemes does not depend on the mean of the prior, the matrix $B$ works as a scalable parameter here.
}

\subsubsection{1-dimensional example}
\begin{figure}
  \centering
  \subfigure{
  \includegraphics[width=0.4\textwidth,height=2.4in]{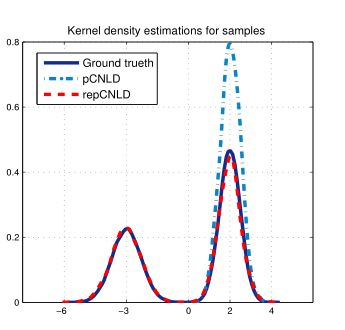}}
  \subfigure{
  \includegraphics[width=0.4\textwidth,height=2.4in]{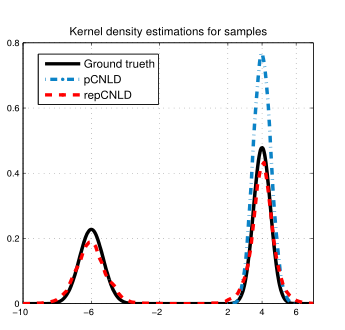}}
  \caption{Kernel density estimations for samples obtained by the pCNLD and repCNLD methods with: ({Left}) $\mu_1=-3, \mu_2=2$ and ({Right}) $\mu_1=-6, \mu_2=4$. The solid line is the true density. }\label{repCNLD1}
\end{figure}

In this experiment, we demonstrate the repCNLD method on a 1-dimensional Gaussian mixture model with ${M}=2$, $\gamma_1=0.4$, $\gamma_2=0.6$, $\theta_1=(-3, 0.7^2)$, $\theta_2=(2, 0.5^2)$, which is the first example in paper \cite{deng2020non}. The time steps are set as $\delta=0.001$, the corresponding values of $\beta$ are calculated as $\beta=0.0447$, and the temperatures $\tau_1$ and $\tau_2$ are set to be 1 and 15, respectively. We also use the pCNLD method to simulate the target distribution as a comparison; the time step and temperature are set to be $\delta_1$ and $\tau_1$. The auxiliary prior is chosen as $\mathcal{N}(0, 3)$ for both methods. We simulate $10^5$ samples from each distribution and the kernel estimations of the samples are displayed in Figure $\ref{repCNLD1}$({Left}). As can be seen, the pCNLD gets stuck in a single mode, while repCNLD successfully yields a close approximation to the ground truth distribution.

We also simulate the distribution with $\mu_1=-6$, $\mu_2=4$, \textcolor{black}{In this example, the distance between the mean of two distributions is significantly larger than in the first example. As a result, it will be harder to capture both modes.} The temperature for the second chain of the repCNLD is set $\tau_2=40$, the auxiliary prior is set $\mathcal{N}(0, 3^2)$, of which the range is wider, and the other parameters are kept the same. The results are shown in Figure $\ref{repCNLD1}$({Right}), the repCNLD can still yield a good approximation, and the pCNLD still fails to escape from the local modes as expected.

\begin{figure}
  \centering
  \subfigure{
  \includegraphics[width=0.4\textwidth,height=2.5in]{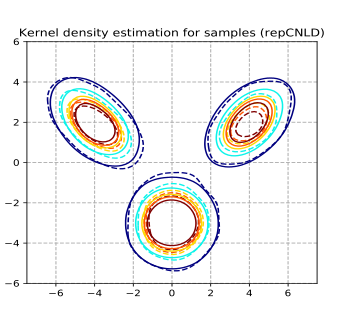}}
  \subfigure{
  \includegraphics[width=0.4\textwidth,height=2.5in]{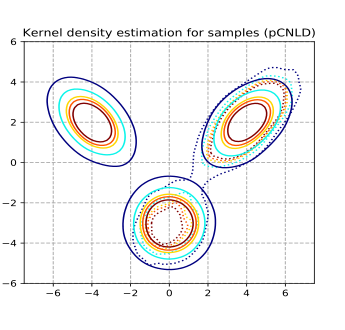}}
  \caption{Density contour plot of ({Left}) samples obtained by the repCNLD method and ({Right}) samples obtained by the pCNLD method. The ground true density is plotted in the solid line in both pictures, the density of samples obtained by the repCNLD method is plotted in the dashed line, while the density of samples obtained by the pCNLD method is shown in dot.}\label{repCNLD2}
\end{figure}

\subsubsection{2-dimensional example}
In this example, we demonstrate the repCNLD method on a 2-dimensional Gaussian mixture model with $M=3$, $\gamma_1=0.3$, $\gamma_2=0.3$, $\gamma_3=0.4$, $m_1=[4, 2]^T$, $m_2=[-4, 2]^T$, $m_3=[0, -3]^T$, and
\[
\begin{gathered}
B_1=
\begin{bmatrix}
1 & 0.6 \\
0.6 & 1
\end{bmatrix},\quad
B_2=
\begin{bmatrix}
1 & -0.6 \\
-0.6 & 1
\end{bmatrix},\quad
B_3=
\begin{bmatrix}
1 & 0 \\
0 & 1
\end{bmatrix}.
\end{gathered}
\]
The time steps are set as $\delta=0.001$, the temperatures $\tau_1$ and $\tau_2$ are set to be 1 and 20 for the repCNLD method, while the time step and temperature are set as $\delta$ and $\tau_1$ for the single chain pCNLD method. The mean and covariance matrix of the Gaussian prior is set as $[0, 0]^T$
and \( B=\left[\begin{array}{cc}
10 & 0 \\
0 & 10
\end{array}\right]\)
. We draw $10^5$ samplers for both the repCNLD and pCNLD methods, the contour plot of the kernel density estimations are illustrated in Figure $\ref{repCNLD2}$. \textcolor{black}{The samples drawn by the repCNLD method provide a good approximation of the target distribution; they explore all modes of the target distribution.} The samples obtained by the pCNLD get trapped in local modes {\color{black}with mean $m_1$ and $m_3$, in addition, few samples are exploring the mode with mean $m_2$.}

\subsection{Initial center identification}

In this section, we consider apply the repCNLD method to nonlinear inverse problems \textcolor{black}{whose underlying forward model is} the single-phase Darcy flow. The flow is described in terms of concentration $u(x,t)$, which is governed by the following parabolic equation
  \begin{equation}
  \label{ex1}
  \frac{\partial u(x,t)}{\partial t}  =\text{div}(\kappa(x)\nabla u(x,t))+f(x,t),\ x\in\Omega,t\in(0,T]\\
  \end{equation}
subject to an appropriate boundary condition and initial condition. The studied physical domain is $\Omega=[0, 1]\times [0, 1]$, the terminal time is set as $T=0.1$.

The initial condition has the form
\[
u_0(x; {\bm \xi})=\frac{q}{2\pi l_{x_1}l_{x_2}}\exp\left(-\frac{(x_1-\xi_1)^2}{2l_{x_1}^2}-\frac{(x_2-\xi_2)^2}{2l_{x_2}^2}\right),
\]
the unknown parameter ${\bm \xi}=(\xi_1, \xi_2)$ is the location of the pollution, the other parameters are prescribed as
\[
l_{x_1}=0.1, \ l_{x_2}=0.2, \ q=1.
\]
The accurate solution of the forward problem is set to be
\[
u(x, t; {\bm \xi})= u_0(x; {\bm\xi})e^{-t},
\]
$\kappa(x)=1$, the inhomogeneous Dirichlet boundary condition and source term are preset to match the solution. Our goal is to infer ${\bm \xi}$ from some measurements of $u(x,t)$ at time $T$. We observe the data at the place $x_o=(0.5, 0.3)$, the measurement error is set as $\sigma_o=0.1$. In this experiment, $u(x_o, T)=M e^{-r^2}e^{-T}$, where $r=1$. This setting guarantees an infinite number of solutions to the inverse problem \cite{lin2021multi}. We set the mean and covariance of the prior as $m= [1/2, 1/2]^T$ and \( B=\left[\begin{array}{cc}
1/4 & 0 \\
0 & 1/4
\end{array}\right]\) to guarantee the domain $\Omega$ be sampled with a certain probability.

In order to show the efficiency and accuracy of the methods, we compare the pCNLD, the repCNLD proceeded with the same forward solver, and the m-repCNLD. We use the standard finite element method to solve the forward models, the backward Euler method is used for the temporal discretization, and the time step is set to be $\Delta t = 0.001$ for all methods. For the m-repCNLD, the solver with $\Delta x = 1/40$ is assigned to the chain with low temperature, while in the high-temperature chain, the forward model is solved with $\Delta x = 1/20$. $\Delta x = 1/40$ is set for both the single chain pCNLD and the repCNLD methods.

The time steps used in discretizing the Langevin diffusion chains are set as $\delta=10^{-5}$, the temperatures are set as $\tau_1=1$, $\tau_2=15$ for both the repCNLD and m-repCNLD methods. For the pCNLD method, the parameters are prescribed as $\delta$ and $\tau_1$. Starting with the same initial value, chains with $3 \times 10^5$ iterations are run. We note that running the chain for $10^5$ iterations, \textcolor{black}{the pCNLD simulation time is about $41462.9$s, while the repCNLD and m-repCNLD methods take $49549.6$s and $47660.7$s to complete all iterations, respectively.} Actually, when the number of iterations is set the same, {\color{black}simulations for the repCNLD or m-repCNLD methods double the simulations for the single chain pCNLD method.} However, the computational cost of the repCNLD method is much less than twice the consumption of the single-chain pCNLD method; this is because the samples from the two chains of the repCNLD or m-repCNLD are updated in parallel. Meanwhile, the parallel technique makes the speeding-up effects of the m-repCNLD not significant. If the parallel toolbox is unavailable, the multi-variance strategy can significantly save the computational cost significantly \cite{lin2021multi}.

\begin{figure}
  \centering
  \includegraphics[width=0.7\textwidth,height=3in]{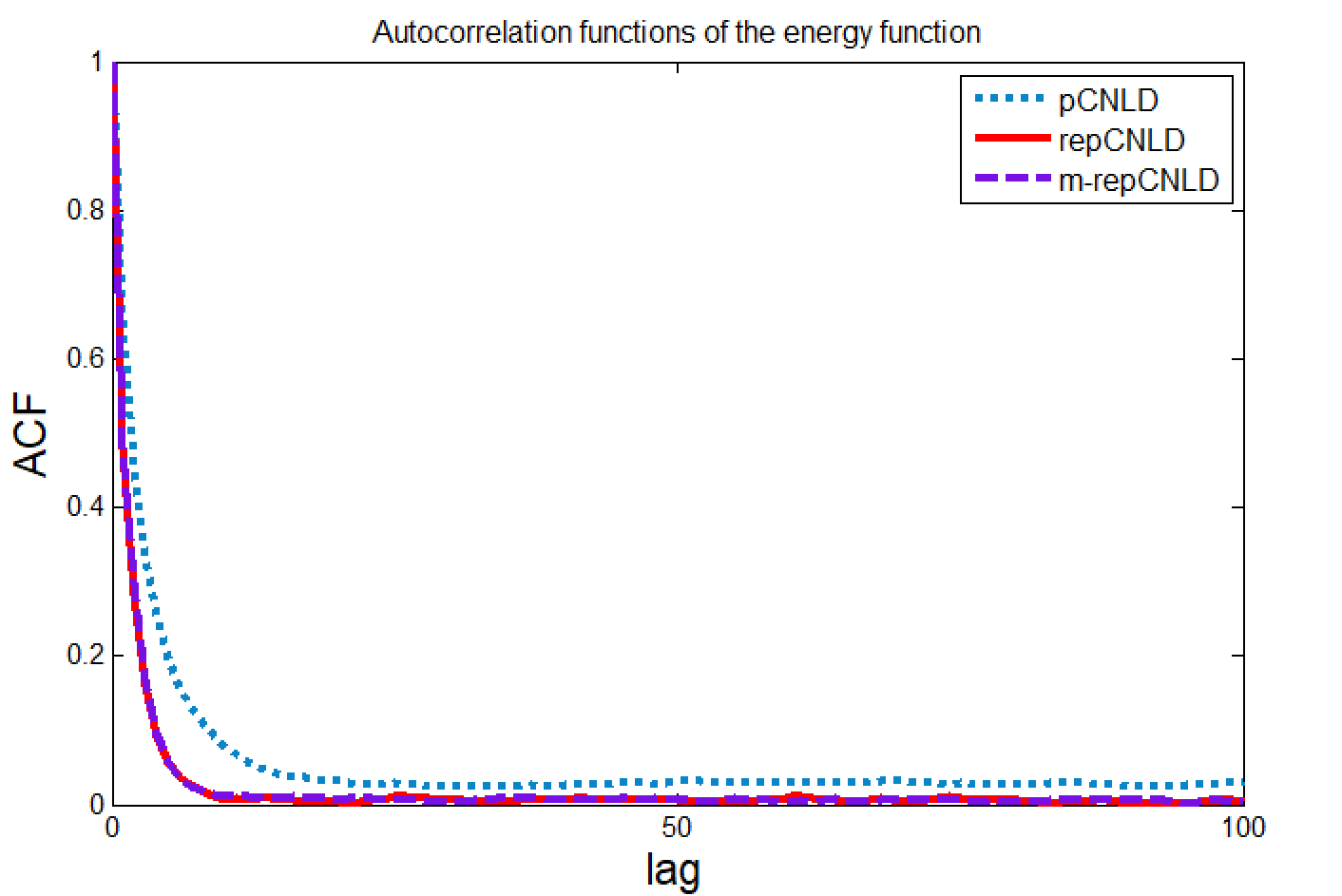}
  \caption{Autocorrelation functions (ACF) of the Onsager-Machlup functional (OMF) or the energy function for different methods. The values of ACF for the replica exchange methods decrease dramatically to 0, and the ACF for the single chain method decreases more slowly compared with the double chain methods.}\label{ACF1}
\end{figure}

\begin{figure}
  \centering
  \includegraphics[width=0.7\textwidth,height=3in]{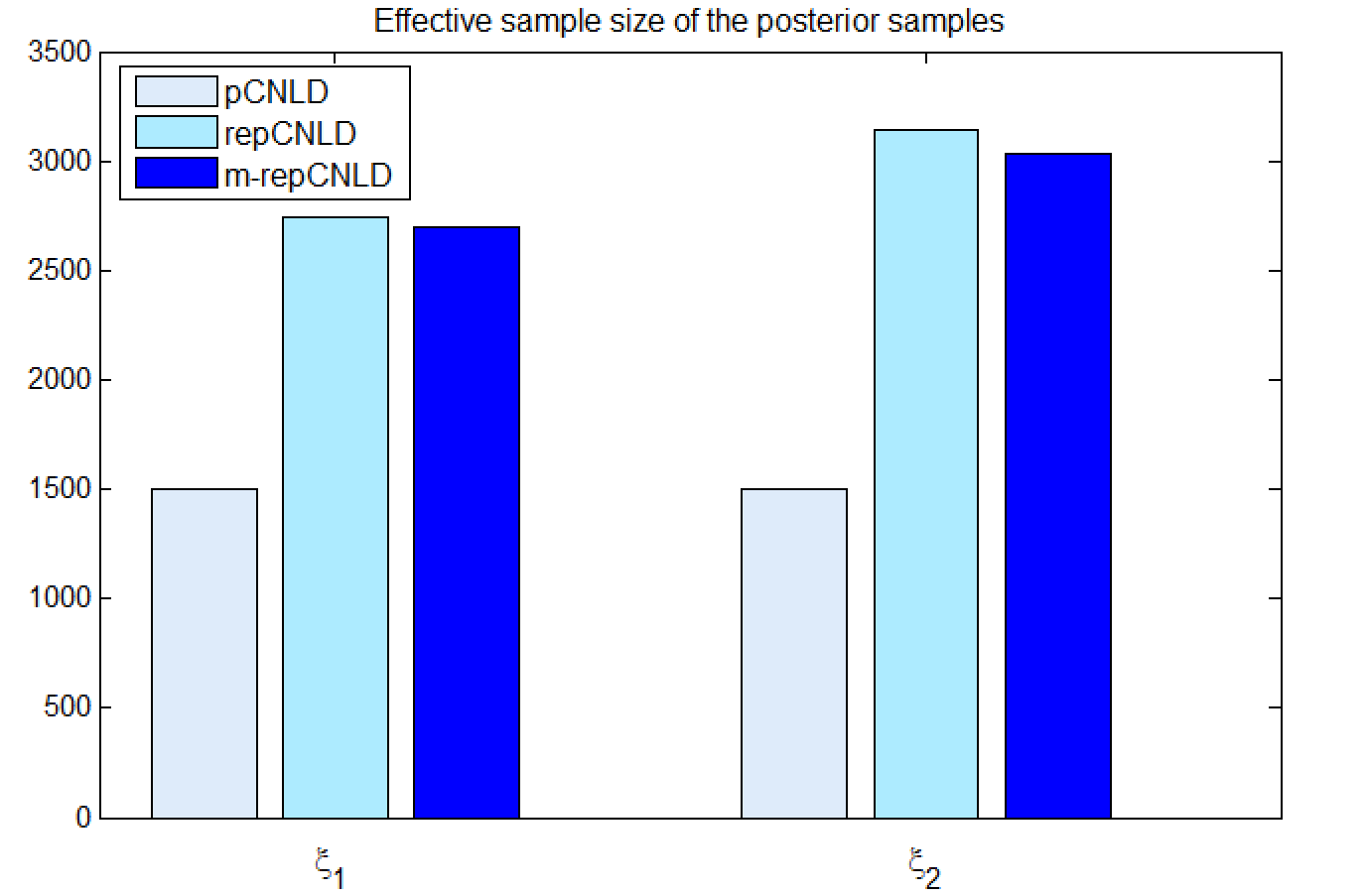}
  \caption{Effective sample size of the posterior samples for different methods. With the same length of Markov chain, the independent samples obtained by the pCNLD method is less than the ones obtained by the replica exchange methods for both $\xi_1$ and $\xi_2$. }\label{ESS}
\end{figure}

\begin{figure}
  \centering
  \includegraphics[width=1\textwidth,height=2.2in]{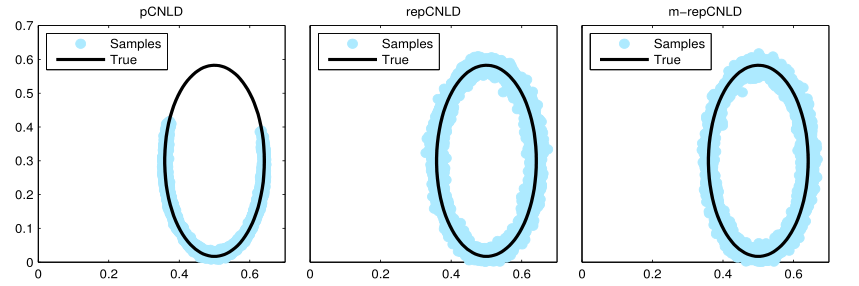}
  \caption{Scatter plot of posterior samples for (Left) pCNLD method, (Middle) repCNLD method, and (Right) m-repCNLD method, respectively. The solid line is the true solution to the inverse problem. }\label{Scatter1}
\end{figure}

We use the autocorrelation functions (ACF) of the Onsager-Machlup functional (OMF)\cite{dashti2013map, li2015note} to compare the efficiency of the three mentioned methods, the OMF here is the corresponding energy function we defined in equation ($\ref{energy}$). The ACF of the OMF for different methods are displayed as functions of the lag in Figure $\ref{ACF1}$. It can be seen from the figure that all of the methods result in low ACF values. The ACF values of repCNLD and m-repCNLD methods decrease much more dramatically than the ACF of the pCNLD method, the values are also much lower. Both the repCNLD and m-repCNLD methods perform better than the pCNLD method.

The effective sample size (ESS) is another scholar to measure the sampling efficiency of MCMC \cite{kass1998markov}. The most common definition of ESS is
\[
\text{ESS}=\frac{N_{max}}{1+2\rho},
\]
where $\rho$ is the integrated autocorrelation time and $N_{max}$ is the total sample size, and it is used to estimate the number of effectively independent samples in the Markov chain. With the same length of the Markov chain, the bigger the ESS is, the more effective the algorithm is. We compute the ESS of the unknown $\xi$ and show them in Figure $\ref{ESS}$, the bar-plot indicates that the replica exchange method generates much more effectively independent draws than the single chain method, while the repCNLD method outperforms the m-repCNLD method in both $\xi_1$ and $\xi_2$.

In Figure $\ref{Scatter1}$, the last $3\times 10^4$ samples of the posterior are plotted for each method, and the solid line is the true solution. We can see that the pCNLD method fails to capture all of the solutions, while the repCNLD and m-repCNLD methods capture all of the solutions to the inverse problem. The second chain of the replica exchange methods successfully helps the particles in the first chain escape from the local modes as expected. In summary, in this inverse problem with infinite solutions, the replica exchange methods are significantly more efficient than the single-chain method.

\subsection{Permeability identification}

In this subsection, we consider the permeability identification of the flow. The Gaussian process is imposed as the prior of the log permeability. Moreover, we adopt the truncated Karhunen Lo\`{e}ve expansion to parameterize the Gaussian field $Y(x, \omega)=\log(\kappa(x; \bm \xi))$ and it can be expressed as
\[
Y(x, \omega)=\bar Y(x) + \sum_{i=1}^n \sqrt{\lambda_i}\varphi(x)\xi_i(\omega),
\]
where $\bar Y(x)$ is the mean of the random field, $\{ \xi_i(\omega)\}_{i=1}^n $ are uncorrelated random variables, and $(\lambda_i, \varphi_i(x))$ are the eigenvalue and eigenfunctions of the eigenvalue problem
\[
\int_\Omega \text{cov}[Y](x_1, x_2)\varphi(x_2)dx_2=\lambda \varphi(x_1),
\]
where $\text{cov}[Y]$ is the covariance function of $Y(x, \omega)$. In this experiment, the mean of the random field is set as $\bar Y(x)=0$, the covariance function of $Y(x, \omega)$ is given by \cite{williams2006gaussian}
\[
\text{cov}[Y](x_1, x_2) = \sigma^2_Y\frac{2^{1-\nu}}{\Gamma(\nu)}\left(\sqrt{2\nu}\frac{|x_1-x_2|}{l}\right)^\nu
B_\nu\left(\sqrt{2\nu}\frac{|x_1-x_2|}{l}\right)
\]
where $\Gamma(\cdot)$ is the Gamma function, and $B_\nu(\cdot)$ is the modified Bessel function. We choose $\sigma^2_Y=1$, $\nu=0.2$, the length scale along the horizontal and vertical direction is set as 1 and 0.5, respectively. $86\%$ of the total prior energy defined by
\[
e(n)=\frac{\sum_{i=1}^n \lambda_i}{\sum_{i=1}^\infty \lambda_i}
\]
is retained, where the eigenvalues are sorted in ascending order, resulting in $n=15$.

\begin{figure}
  \centering
  \subfigure{
  \includegraphics[width=0.45\textwidth,height=2.3in]{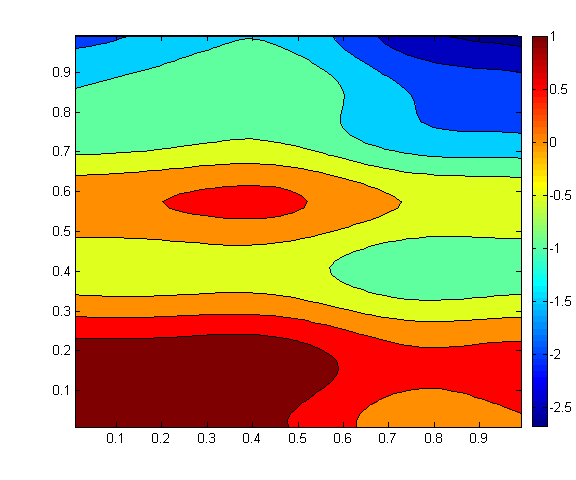}}
  \subfigure{
  \includegraphics[width=0.45\textwidth,height=2.3in]{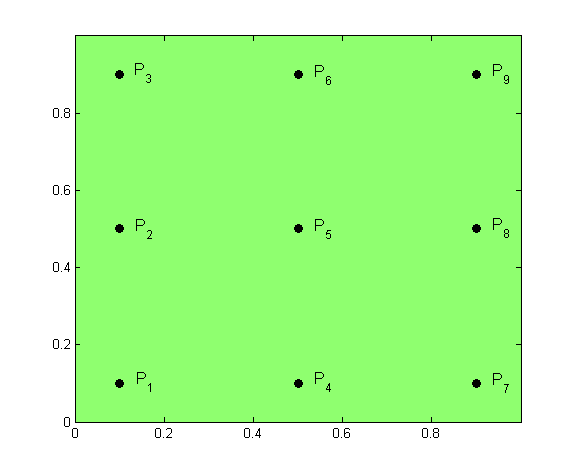}}
  \caption{({Left}) The true log permeability, ({Right}) well configuration.}\label{ref-locs}
\end{figure}

The force term in the equation $(\ref{ex1})$ is set as
\[
f(x, t)=\sum_{i=1}^{N_w}\frac{q_i}{2\pi l_w^2}\exp\bigg\{-\frac{\|x-x_i\|^2}{2l_w^2}\bigg\},
\]
where $\{q_i\}_{i=1}^{N_w}$ and $\{x_i\}_{i=1}^{N_w}$ are the prescribed production rates and well locations \cite{chen2006computational}. The production rate is set as 1.5 and $l_w=0.05$ in the force term. The true coefficient $\kappa(x; \bm \xi)$ is a realization from the prior, the true log permeability and the well locations are shown in Figure $\ref{ref-locs}$. In addition, we consider the initial condition $u(x, 0) = 4$ and no-flow boundary conditions for this problem.

The synthetic data is generated by using FEM on a $120\times 120$ grid and 100-time steps. The observed data is the pressure collected at the well locations, the measurement time used are $t_1=0.01$, $t_n=0.02n$ ($n=2, \cdots, 5$), the observed noise is $\sigma_o=0.1$.   The forward model is solved with $\Delta x = 1/60$, and the time interval $[0, T]$ is discretized with 50 steps for the repCNLD and pCNLD methods. For the m-repCNLD method, the forward model is solved with $\Delta t=0.002$ for the first chain, while the forward model simulated in the second chain is run by $\Delta t=0.005$.

\begin{figure}
  \centering
  \includegraphics[width=0.7\textwidth,height=3in]{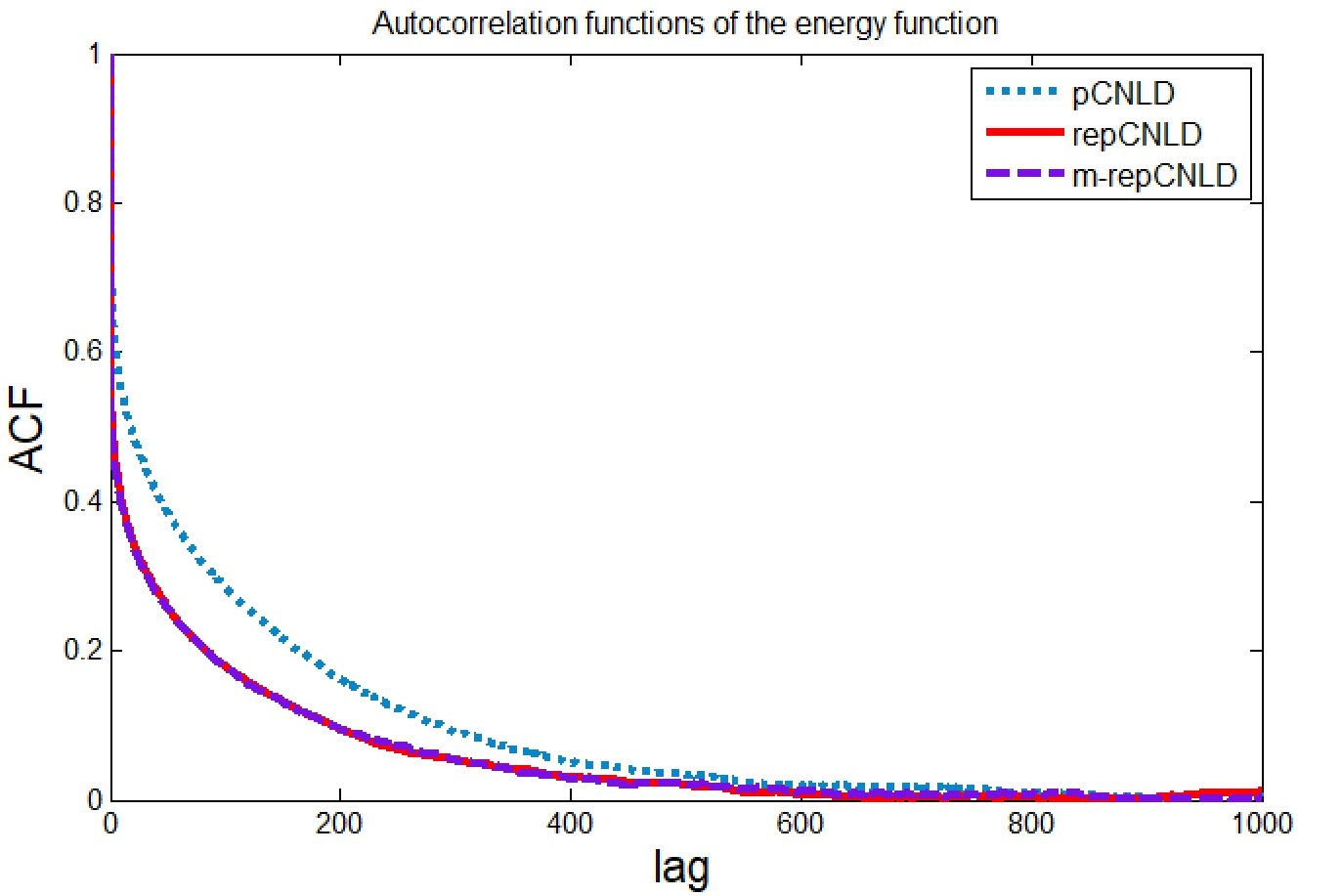}
  \caption{Autocorrelation functions (ACF) of the Onsager-Machlup functional (OMF) or the energy function for different methods. The values of ACF for the replica exchange methods are lower than the ones of the pCNLD method, the chains obtained by the replica exchange methods have a higher mixing rate than the ones obtained by the pCNLD method. }\label{ACF2}
\end{figure}

The time steps involved in the Langevin dynamics are set as $\delta=2\times10^{-3}$, the temperatures are set as $\tau_1=1$, $\tau_2=1.5$ for both the repCNLD and m-repCNLD methods. The parameters are prescribed as $\delta_1$ and $\tau_1$ for the pCNLD method. The starting point is set as the mean of the prior. The computational cost of $10^5$ iterations for the pCNLD, repCNLD and m-repCNLD methods is $57263.6$s, $58535.2$s, $57212.7$s, respectively. This is also due to the usage of the adjoint method and parallel technique. $5\times 10^5$ samples are drawn for each method, the ACF of the OMF is plotted against the lag in Figure $\ref{ACF2}$, the ACF values of the repCNLD and m-repCNLD methods are lower than the pCNLD method, which means that they have higher mixing rate than the pCNLD method.

\begin{figure}
  \centering
  \includegraphics[width=0.6\textwidth,height=3.8in]{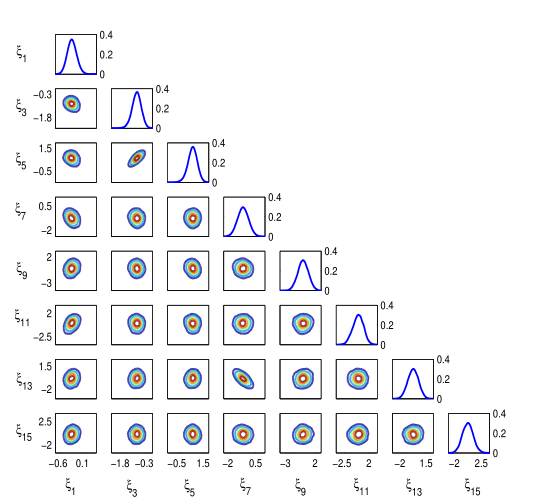}
  \caption{The 1-D and 2-D marginal posterior probability distributions for part of the parameters for the pCNLD method. The uncertainty of the modes $[\xi_1, \xi_3, \xi_5]$ is reduced by the observed data, while for the modes that are associated with small eigenvalues, the data is not informative.}\label{kde-pCNLD}
\end{figure}

\begin{figure}
  \centering
  \includegraphics[width=0.6\textwidth,height=3.8in]{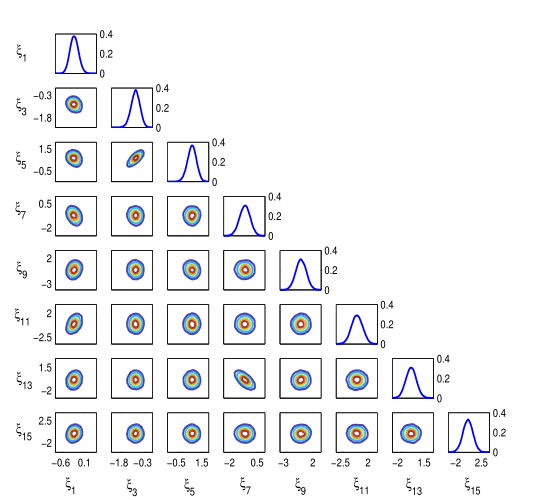}
  \caption{The 1-D and 2-D marginal posterior probability distributions for part of the parameters for the repCNLD method.}\label{kde-repCNLD}
\end{figure}

\begin{figure}
  \centering
  \includegraphics[width=0.6\textwidth,height=3.8in]{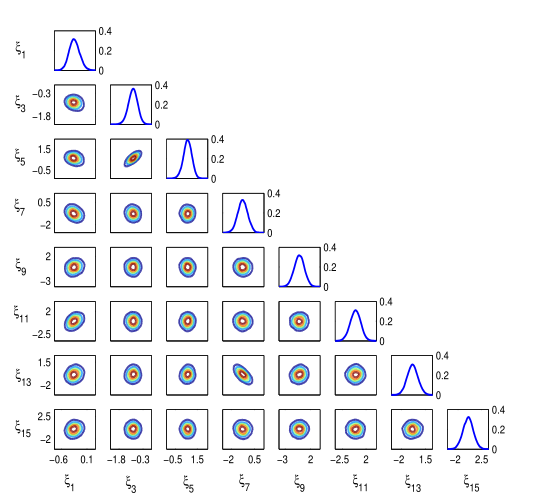}
  \caption{The 1-D and 2-D marginal posterior probability distributions for part of the parameters for the m-repCNLD method.}\label{kde-repCNLDm}
\end{figure}

\begin{figure}
  \centering
  \subfigure{
  \includegraphics[width=1\textwidth,height=2in]{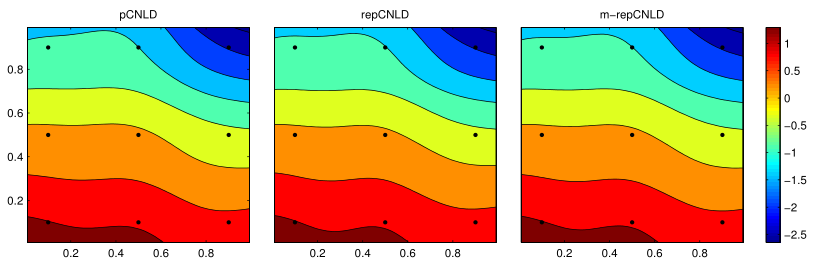}}
  \subfigure{
  \includegraphics[width=1\textwidth,height=2in]{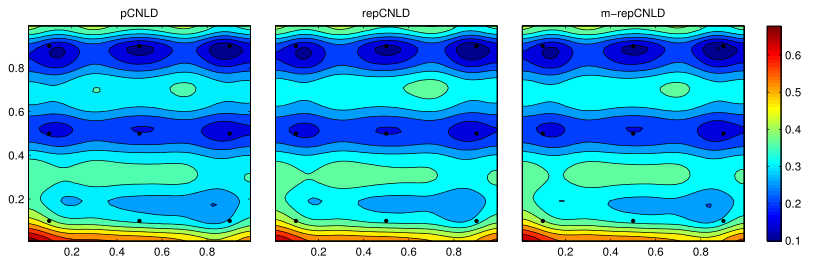}}
  \subfigure{
  \includegraphics[width=1\textwidth,height=2in]{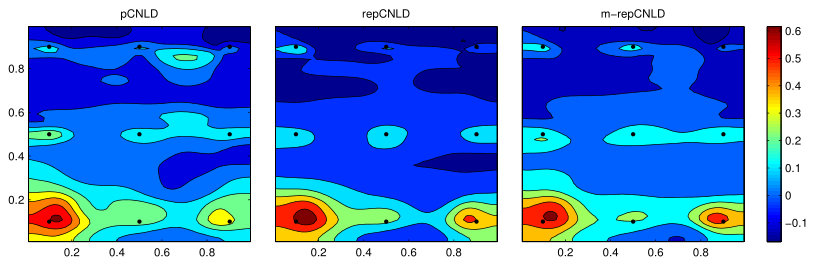}}
  \caption{The estimated posterior (Upper) mean, (Middle) standard deviation and (Bottom) skewness of the log $\kappa(x; \xi)$ for different methods. The spatial distributions of estimates are basically the same. The skewness values around the wells $P_1$ and $P_7$ are larger than the ones around the other wells, since the data provided by the two wells is less informative. }\label{pos}
\end{figure}

We visualize the 1-dimensional (1-D) and 2-dimensional (2-D) marginal posterior probability distributions for part of the parameters for the three methods. Based on the shape of their 2-D marginals, some correlations are apparent between several modes, such as $\xi_1$ and $\xi_7$, $\xi_3$ and $\xi_5$, $\xi_7$ and $\xi_{13}$. The ranges of the 1-D marginal posterior for parameters $\xi_1$, $\xi_3$ and $\xi_5$ are much narrower than the ones of the high-indexed modes, which indicates that the infinite data is not informative about the small scales of the unknown log permeability. It can be observed from Figure $\ref{pos}$, the posterior mean and standard deviation for all the methods have common spatial features. There is more uncertainty around the wells $P_1$ and $P_4$, and less uncertainty around the other wells. The non-zero values of the posterior skewness indicate that the posterior is non-Gaussian; this may be contributed by the non-linearity of the forward map and the insufficient data information.

\section{Conclusions}
In this work, we propose the repCNLD MC algorithm to explore target or posterior distributions; two preconditioned Langevin diffusion chains with different temperatures are simulated for the sampling process. The diffusion with low-temperature aims to exploit the geometry locally, while the diffusion with high temperature works to explore the domain globally. They are allowed to swap with the specified probability. The swapping does not change the invariant distribution of the diffusion chains but accelerates the convergence of the target distribution. We use the Crank-Nicolson scheme to discretize the  preconditioned Langevin diffusion chains; the covariance of the Gaussian prior is cooperated with the temperature parameters to scale the noise of the diffusion together. We analyze the discretization error of using the pCN scheme. The proposed method performs well for multimodal distribution taking advantage of the replica exchange method. In order to save the computational cost, the forward model applied in the second chain is solved with the low-fidelity method; we derive an unbiased estimator for the swapping exchange rate by tracing to the forward model errors. Meanwhile, we provide the framework of a discrete adjoint method for the calculation of the gradient. We demonstrate the efficiency of the proposed method through some numerical examples. For the nonlinear inverse problems, we use the ACF values of OMF (energy function) values to assess the MCMC chains' efficiency. The replica exchange methods outperform the single-chain methods with the same starting points.

We use the adjoint method to calculate the derivative information of the log-likelihood function, which saves the computational cost. However, it would be unattainable for problems with black-box models. The covariance matrix of the prior works is a preconditioned matrix here; it is preset and fixed; some adaptive methods may be used to improve the efficiency of our proposed method to deal with complex problems where the scales of the objective function are different in different modes. Further investigation of these issues is worth pursuing in the future.

\section*{Acknowledgments}
NO acknowledges the support of Hunan Provincial Natural Science Foundation of China 2021JJ40557 and Chinese NSF 11901060.
GL gratefully acknowledges the support of the National Science Foundation (DMS-1555072, DMS-2053746, and DMS-2134209), and Brookhaven National Laboratory Subcontract 382247, and U.S. Department of Energy (DOE) Office of Science Advanced Scientific Computing Research program DE-SC0021142 and DE-SC0023161.
\smallskip
\bigskip

\bibliographystyle{plain}
\bibliography{pCN-replica}

\section*{Appendix}

\subsection*{Proof of theorem \ref{invariance}}
\begin{proof}
The dynamics in equations $(\ref{cns-U})$ and $(\ref{swap})$ defines a Markov jump process, following \cite{oksendal2013stochastic} and \cite{chen2020accelerating}, the infinitesimal generator has the form
\[
\begin{array}{l}
\mathcal{L}^c\left(f(x, y)\right)=\underbrace{-[L\nabla_{x}U(x)] \cdot  \nabla_{x}f(x, y)+\tau_{1} L : \nabla_{x}\nabla_{x} f(x, y)}_{\mathcal{L}^c_1\left(f(x, y)\right)}\\
\underbrace{-[L\nabla_{y}U(y)] \cdot  \nabla_{y}f(x, y)+\tau_2 L : \nabla_{y}\nabla_{y} f(x, y)}_{\mathcal{L}^c_2\left(f(x, y)\right)}+\underbrace{ cs(x, y)\left(f(y, x)-f(x, y)  \right)}_{\mathcal{L}^c_s\left(f(x, y)\right)},
\end{array}
\]
As demonstrated in the proof of the Lemma 3.2 in paper \cite{chen2020accelerating}, we only need to show that
\begin{equation}\label{detailbalance}
\int \int g(x, y) \mathcal{L}^{c}\left(f(x, y)\right) \pi(x, y){\rm d}x{\rm d}y=\int \int f(x, y) \mathcal{L}^{c}\left(g(x, y)\right) \pi(x, y){\rm d}x{\rm d}y
\end{equation}
for all $f, g \in \mathcal{D}(\mathcal{L}^c)$, where the density $\pi(x, y)$ is defined as $(\ref{invariant})$, to prove the invariance and
reversibility. Since the term arises from the swap $\mathcal{L}_s^{c}$ is the same with the one in paper \cite{chen2020accelerating}, we directly have
\[
\int \int g(x, y) \mathcal{L}_s^{c}\left(f(x, y)\right) \pi(x, y){\rm d}x{\rm d}y=\int \int f(x, y) \mathcal{L}_s^{c}\left(g(x, y)\right) \pi(x, y){\rm d}x{\rm d}y.
\]
For the term $\mathcal{L}_1^{c}$, by integration by parts and the fact that $f$ and $g$ have compact support, we have that for all fixed $y \in \mathbb{R}^n$,
\setlength{\arraycolsep}{1.4pt}
\[
\begin{array}{l}
\mathlarger{\int} \tau_1g(x, y)\left[L: \nabla_{x}\nabla_{x} f(x, y)\right]\exp(\frac{U(x)}{\tau_1}){\rm d}x \\
=-\tau_1 \mathlarger{\int} L : \left(\nabla_{x} f(x, y) \left(\nabla_{x}\left[g(x, y)\exp(\frac{U(x)}{\tau_1})\right]\right)^T\right){\rm d}x\\
= -\mathlarger{\int} \tau_1\exp(\frac{U(x)}{\tau_1})\left[L: \left(\nabla_{x}f(x, y)\left(\nabla_{x}g(x, y)\right)^T\right)\right]{\rm d}x\\
   \quad +\mathlarger{\int} g(x, y)\exp(\frac{U(x)}{\tau_1})\left[L: \left(\nabla_{x}f(x, y)\left(\nabla_{x}U(x)\right)^T\right)\right]{\rm d}x, \\
\end{array}
\]
where the symbol $T$ means the transposition of a matrix. Note that the term \([L\nabla_{x}U(x)] \cdot  \nabla_{x}f(x, y) \) can be rewritten as \(L : \left(\nabla_{x}f(x, y)  \left(\nabla_{x} U(x) \right)^T \right) \), we then have
\[
\int g(x, y) \mathcal{L}_1^{c}\left(f(x, y)\right) \pi(x, y){\rm d}x=-\int \tau_1\left[L: \left(\nabla_{x}f(x, y)\left(\nabla_{x}g(x, y)\right)^T\right)\right]\pi(x, y){\rm d}x,
\]
hence,
\[
\int\int g(x, y) \mathcal{L}_1^{c}\left(f(x, y)\right) \pi(x, y){\rm d}x{\rm d}y=-\int\int \tau_1\left[L: \left(\nabla_{x}f(x, y)\left(\nabla_{x}g(x, y)\right)^T\right)\right]\pi(x, y){\rm d}x{\rm d}y.
\]
By switching the positions of $f$ and $g$ in the above equation, we have
\[
\int\int f(x, y) \mathcal{L}_1^{c}\left(g(x, y)\right) \pi(x, y){\rm d}x{\rm d}y=-\int\int \tau_1\left[L: \left(\nabla_{x}g(x, y)\left(\nabla_{x}f(x, y)\right)^T\right)\right]\pi(x, y){\rm d}x{\rm d}y,
\]
since the matrix $L$ is symmetric, then
\[
L: \left(\nabla_{x}f(x, y)\left(\nabla_{x}g(x, y)\right)^T \right)=L: \left(\nabla_{x}g(x, y)\left(\nabla_{x}f(x, y)\right)^T \right),
\]
we finally have
\begin{equation}\label{ig1}
\int \int g(x, y) \mathcal{L}_1^{c}\left(f(x, y)\right) \pi(x, y){\rm d}x{\rm d}y=\int \int f(x, y) \mathcal{L}_1^{c}\left(g(x, y)\right) \pi(x, y){\rm d}x{\rm d}y.
\end{equation}
By the same derivation, the equation $(\ref{ig1})$ also holds for $\mathcal{L}_2^{c}$. In summary, the equation $(\ref{detailbalance})$ is proved to be established for the density $\pi$ defined in $(\ref{invariant})$. Thus, $\{\bm{\xi}_t\}_{t\geq 0} $ is a reversible Markov process, and its invariant distribution has the density $\pi(\bm{\xi}^{(1)}, \bm{\xi}^{(2)})\propto \exp\{-\frac{U(\bm{\xi}^{(1)})}{\tau_1}-\frac{U(\bm{\xi}^{(2)})}{\tau_2}\}$.
\end{proof}

\subsection*{Proof of theorem \ref{acceleration}}
\begin{proof}
According to the property of $\nabla$, we have
\begin{eqnarray*}
  \nabla_{z}f^2 &=& 2f \nabla_{z}f,\\
  \nabla_{z}\nabla_{z}f^2 &=& 2\left(\nabla_{z}f\left(\nabla_{z}f \right)^T\right)+2f\left(\nabla_{z}\nabla_{z}f\right)
\end{eqnarray*}
for $z=x, y$. Then we have
\[
\begin{array}{l}
\Gamma^{c}\left(f\left(x, y\right)\right)=\frac{1}{2}\mathcal{L}^{c}\left(f^{2}\left(x, y\right)\right)-f\left(x, y\right)\mathcal{L}^{c}\left(f\left(x, y\right)\right)\\
=\tau_1\left[L\nabla_{x}f(x, y)\right]\cdot\nabla_{x}f(x, y)+\tau_2\left[L\nabla_{y}f(x, y)\right]\cdot\nabla_{y}f(x, y)\\
\quad+\frac{c}{2} s(x, y)\left(f(y, x)-f(x, y)\right)^2,\\
\end{array}
\]
since the matrix $L$ is symmetric positive definite, the first two terms of the above equation is consequently non-negative. Hence, the Dirichlet form $\mathcal{E}^c$ can be given by
\[
\begin{aligned}
\mathcal{E}^{c}(f)=\int &\left(\tau_{1}\left\|\sqrt{L}\nabla_{x} f(x, y)\right\|^{2}+\tau_{2}\left\|\sqrt{L}\nabla_{y} f(x, y)\right\|^{2}\right) \mathrm{d} \mu\left(x, y\right) \\
&+\int \frac{c}{2} s\left(x, y\right)\left(f(y, x)-f(x, y)\right)^{2} \mathrm{d} \mu\left(x, y\right)
\end{aligned}
\]
\end{proof}

\subsection*{Proof of theorem \ref{discrete-error}}
\begin{proof}
The integral form of equation $(\ref{continue})$ is
\[
\bm{\xi}_t =\bm{\xi}_0-\int_0^t H\nabla D(\bm{\xi}_s){\rm d}s+\int_0^t \sqrt{H}\Sigma_s\mathrm{d}\bm{W}_s,
\]
compare it with the equation $(\ref{continue-interpolation})$, we have
\[
\begin{aligned}
\bm{\xi}_t-\bm{\xi}_t^\delta&=-H\int_0^t\left(\nabla D(\bm{\xi}_t)-\nabla D\left({\bm \xi}^{\delta}_{\lfloor s/{\delta} \rfloor \delta}\right) \right)\mathrm{d}s-(1-\eta)H\int_0^t \nabla D\left({\bm \xi}^{\delta}_{\lfloor s/{\delta} \rfloor \delta}\right){\rm d}s\\
&+\sqrt{H}\int_0^t \left(\Sigma_s-\Sigma^{\delta}_{\lfloor s/{\delta} \rfloor \delta}\right) {\rm d}\bm{W}_s
+(1-\eta)\sqrt{H}\int_0^t \Sigma^{\delta}_{\lfloor s/{\delta} \rfloor \delta}{\rm d}\bm{W}_s.
\end{aligned}
\]
According to the compatibility of matrix norms, we apply the Cauchy-Schwartz inequality and take expectation, then
\begin{equation}\label{total}
\begin{aligned}
{\bb E}\left[ \|\bm{\xi}_t-\bm{\xi}_t^\delta\|^2 \right]&\leq 4\|H\|_F^2\underbrace{{\bb E}\left[\bigg\|\int_0^t\left( \nabla D(\bm{\xi}_s)-\nabla D\left({\bm \xi}^{\delta}_{\lfloor s/{\delta} \rfloor \delta}\right)\right)\mathrm{d}s\bigg\|^2\right]}_{{\cal I}_1}\\
&+4tr(H)\underbrace{{\bb E}\left[\bigg\|\int_0^t \left(\Sigma_s-\Sigma^{\delta}_{\lfloor s/{\delta} \rfloor \delta}\right) {\rm d}\bm{W}_s\bigg\|^2\right]}_{{\cal I}_2}\\
&+4(1-\eta)^2\|H\|_F^2\underbrace{{\bb E}\left[\bigg\|\int_0^t \nabla D\left({\bm \xi}^{\delta}_{\lfloor s/{\delta} \rfloor \delta}\right){\rm d}s \bigg\|^2\right]}_{{\cal J}_1}\\
&+4(1-\eta)^2tr(H)\underbrace{{\bb E}\left[\bigg\|\int_0^t \Sigma^{\delta}_{\lfloor s/{\delta} \rfloor \delta}{\rm d}\bm{W}_s\bigg\|^2\right]}_{{\cal J}_2}.
\end{aligned}
\end{equation}
where $tr(H)$ is the trace of the matrix $H$, the equality $\|\sqrt{H}\|_F^2=tr(H)$ establishes.

For the term ${\cal J}_1$, by Cauchy-Schwarz inequality and the smooth assumption, we have
\[
\begin{aligned}
{\cal J}_1&={\bb E}\left[ \bigg\|\int_0^t\left(\nabla D\left({\bm \xi}^{\delta}_{\lfloor s/{\delta} \rfloor \delta}\right)- \nabla D(\bm{\xi}^*)\right)\mathrm{d}s\bigg\|^2\right]
\leq t{\bb E}\left[ \int_0^t\bigg\|\left(\nabla D\left({\bm \xi}^{\delta}_{\lfloor s/{\delta} \rfloor \delta}\right)- \nabla D(\bm{\xi}^*)\right)\bigg\|^2\mathrm{d}s\right]\\
&\leq tK^2{\bb E}\left[ \int_0^t\left(\|{\bm \xi}^{\delta}_{\lfloor s/{\delta} \rfloor \delta}-\bm{\xi}^*\|^2\right)\mathrm{d}s\right]
\leq 2tK^2{\bb E}\left[ \int_0^t\left(\|{\bm \xi}^{\delta}_{\lfloor s/{\delta} \rfloor \delta}\|^2+\|\bm{\xi}^*\|^2\right)\mathrm{d}s\right]\\
&\leq 2tK^2 {\bb E}\left[ \sum_{k=0}^{\lfloor t/{\delta} \rfloor} \int_{k \delta}^{(k+1)\delta} \left(\|{\bm \xi}^{\delta}_{k \delta}\|^2+\|\bm{\xi}^*\|^2\right)\mathrm{d}s\right]\\
&\leq 2tK^2(1+t/\delta)\delta\left(\sup_{k\geq 0} {\bb E}\left[\|{\bm \xi}^{\delta}_{k \delta}\|^2\right]+\|\bm{\xi}^*\|^2\right)
\end{aligned}
\]
where $\bm{\xi}^*$ is a stationary point of $D(\cdot)$, i.e., $\nabla D(\bm{\xi}^*)=0$.  Analogue to the Lemma C.2 in paper \cite{chen2020accelerating}, when \( 0<\frac{2\delta}{2+\delta}<\frac{\alpha_1}{\|B\|_F^2K^2}\), we have
\begin{equation}\label{sub-bound}
\sup_{k\geq 0} {\bb E}\left[\|{\bm \xi}^{\delta}_{k \delta}\|^2\right] \leq \epsilon_1 (n, \tau_2, K, \alpha_1, \alpha_2, B)
\end{equation}
where the constant $\epsilon_1(n, \tau_2, K, \alpha_1, \alpha_2, B)$ depends on the dimension $n$, the temperature parameter $\tau_2$, the smoothness constant $K$, the constants $(\alpha_1, \alpha_2)$ in equation $(\ref{B-dissipative})$ and the covariance matrix $B$. Hence, we have
\begin{equation}\label{J1}
{\cal J}_1\leq 2tK^2(1+t/\delta)\delta\epsilon_2 (n, \tau_2, K, \alpha_1, \alpha_2, B).
\end{equation}

For the term ${\cal J}_2$, according to the It$\hat{o}$ isometry, we have
\[
\begin{aligned}
{\cal J}_2 &=\sum_{j=1}^{2n}\int_0^t  {\bb E}\left[\left(\Sigma^{\delta}_{\lfloor s/{\delta} \rfloor \delta}(j)\right)^2\right] {\rm d}s\\
&\leq \sum_{j=1}^{2n} \sum_{k=0}^{\lfloor t/\delta \rfloor} \int_{k\delta}^{(k+1)\delta} {\bb E}\left[\left(\Sigma^{\delta}_{k \delta}(j)\right)^2\right] {\rm d}s   \\
\end{aligned}
\]
where $\Sigma^{\delta}_{k \delta}(j)$ is the $j$-th diagonal entry of $\Sigma^{\delta}_{k \delta}$. Recall that the values of $\left(\Sigma^{\delta}_{k \delta}(j)\right)^2$ would be $2\tau_1$ or $2\tau_2$, suppose the probability of being $2\tau_2$ is $p\in [0, 1]$, then the probability of being $2\tau_1$ is $1-p$, then
\[
{\bb E}\left[\left(\Sigma^{\delta}_{k \delta}(j)\right)^2\right]=2p(\tau_2-\tau_1)+2\tau_1\leq 2\tau_2,
\]
hence we have
\begin{equation}\label{J2}
{\cal J}_2\leq 4(1+t/\delta)\delta\tau_2 n.
\end{equation}

Similar to the analysis in paper \cite{chen2020accelerating}, we have
\[
\begin{aligned}
{\cal I}_1 &\leq t {\bb E}\left[ \int_0^t \bigg\|\nabla D(\bm{\xi}_s)-\nabla D\left({\bm \xi}^{\delta}_{\lfloor s/{\delta} \rfloor \delta}\right)\bigg\|^2 \mathrm{d}s\right]
\leq tK^2{\bb E}\left[ \int_0^t \|\bm{\xi}_s-{\bm \xi}^{\delta}_{\lfloor s/{\delta} \rfloor \delta}\|^2\mathrm{d}s\right]\\
&\leq 2K^2t\left( {\bb E}\left[ \int_0^t \|\bm{\xi}_s-{\bm \xi}^{\delta}_s\|^2\mathrm{d}s\right]+ {\bb E}\left[ \int_0^t \|\bm{\xi}^{\delta}_s-{\bm \xi}^{\delta}_{\lfloor s/{\delta} \rfloor \delta}\|^2\mathrm{d}s\right]  \right)\\
\end{aligned}
\]
where
\[
{\bb E}\left[ \int_0^t \|\bm{\xi}^{\delta}_s-{\bm \xi}^{\delta}_{\lfloor s/{\delta} \rfloor \delta}\|^2\mathrm{d}s\right]\leq(1+t/\delta)\left(4K^2\|H\|_F^2\delta^3\left(\sup_{k\geq 0} {\bb E}\left[\|{\bm \xi}^{\delta}_{k \delta}\|^2\right]+\|\bm{\xi}^*\|^2\right)+4tr(H)n\tau_2\delta^2 \right),
\]
hence we have
\begin{equation}\label{I1}
{\cal I}_1\leq 2K^2t\left( {\bb E}\left[ \int_0^t \|\bm{\xi}_s-{\bm \xi}^{\delta}_s\|^2\mathrm{d}s\right]+\epsilon_3(n,\tau_2,K,\alpha_1,\alpha_2,B)\delta \right).
\end{equation}

Since the term ${\cal I}_2$ is the same with the one in paper \cite{chen2020accelerating}, we use the result directly and have
\begin{equation}\label{I2}
{\cal I}_2\leq \epsilon_4(n,\tau_1,\tau_2,c)t \delta.
\end{equation}

By substituting the equations $(\ref{J1})$, $(\ref{J2})$, $(\ref{I1})$ and $(\ref{I2})$ into $(\ref{total})$, we have
\[
\begin{aligned}
{\bb E}\left[ \|\bm{\xi}_t-\bm{\xi}_t^\delta\|^2 \right]&\leq 8\|B\|_F^2\cdot 2K^2t\left( {\bb E}\left[ \int_0^t \|\bm{\xi}_s-{\bm \xi}^{\delta}_s\|^2\mathrm{d}s\right]+\epsilon_3(n,\tau_2,K,\alpha_1,\alpha_2,B)\delta \right)\\
&+8(1-\eta)^2(1+t/\delta)\|B\|_F^2\cdot 2tK^2\epsilon_2 (n, \tau_2, K, \alpha_1, \alpha_2, B)\delta \\
&+8(1-\eta)^2(1+t/\delta)tr(B)\cdot 4\tau_2 n\delta +8tr(B)\cdot\epsilon_4(n,\tau_1,\tau_2,c)t \delta,
\end{aligned}
\]
recall that $\eta=2/(2+\delta)$, $\delta \in (0, 2)$, then for $t\in [0, T]$, we have
\[
(1-\eta)^2(1+t/\delta)=\frac{\delta^2+t\delta}{(\delta+2)^2}\leq (\frac{1}{4}+\frac{T}{8}).
\]
Finally, we have
\[
\begin{aligned}
{\bb E}\left[ \|\bm{\xi}_t-\bm{\xi}_t^\delta\|^2 \right]&\leq 8\|B\|_F^2\cdot 2K^2t\left( {\bb E}\left[ \int_0^t \|\bm{\xi}_s-{\bm \xi}^{\delta}_s\|^2\mathrm{d}s\right]+\epsilon_3(n,\tau_2,K,\alpha_1,\alpha_2,B)\delta \right)\\
&+\epsilon_6(n,\tau_1,\tau_2,c,K,\alpha_1,\alpha_2,T,B)\delta.
\end{aligned}
\]
Hence, by applying the Gr$\ddot{o}$nwall's equality, there exists a constant $\epsilon(n,\tau_1,\tau_2,c,K,\alpha_1,\alpha_2,T,B)$ such that for $t\in [0, T]$,
\[
{\bb E}\left[ \|\bm{\xi}_t-\bm{\xi}_t^\delta\|^2 \right]\leq \epsilon(n,\tau_1,\tau_2,c,K,\alpha_1,\alpha_2,T,B)\delta.
\]
\end{proof}
\end{document}